\documentclass[preprint,review,fleqn,sort&compress,3p]{elsarticle}%preprint,review,fleqn,sort&compress,3p
\usepackage{mathrsfs}
\usepackage{amsfonts}
\usepackage{subfigure}
\usepackage{caption}
\usepackage{bm}
\usepackage{mathtools}
\usepackage{pifont}
\usepackage{float}
\usepackage{makecell}
\usepackage{natbib}
\usepackage{geometry}
\usepackage{color,xcolor}
\usepackage{url}
\usepackage{booktabs}
\usepackage{multirow}  
\usepackage{graphics}
\usepackage{amssymb}
\usepackage{amsthm}
\usepackage{amsmath}
\usepackage{arydshln}%虚线表

\newtheorem{lemma}{Lemma}[section]
\newtheorem{proposition}{Proposition}[section]
\newtheorem{example}{Example}
\newtheorem{theorem}{Theorem}[section]
\newtheorem{assumption}{Assumption}
\newtheorem{remark}{Remark}
\journal{Applied Mathematics and Computation}
\begin{document}
\begin{frontmatter}
\title{Double fast algorithm for solving time-space fractional diffusion problems with spectral fractional Laplacian}
\author[1]{Yi Yang}%\snm{}
\ead{eyannyee@163.com}
\author[1]{Jin Huang\corref{cor1}}
\cortext[cor1]{Corresponding author}
\ead{huangjin12345@163.com, huangj@uestc.edu.cn}
\address[1]{School of Mathematical Sciences, University of Electronic Science and Technology of China, Chengdu, 611731, China}
\begin{abstract}
	This paper presents an efficient and concise double fast algorithm to solve high dimensional time-space fractional diffusion problems with spectral fractional Laplacian. We first establish semi-discrete scheme of time-space fractional diffusion equation, which uses linear finite element or fourth-order compact difference method combining with matrix transfer technique to approximate spectral fractional Laplacian. Then we introduce a fast time-stepping L1 scheme for time discretization. The proposed scheme can exactly evaluate fractional power of matrix and perform matrix-vector multiplication at per time level by using discrete sine transform, which doesn't need to resort to any iteration method and can significantly reduce computation cost and memory requirement. Further, we address stability and convergence analyses of full discrete scheme based on fast time-stepping L1 scheme on graded time mesh. Our analysis shows that the choice of graded mesh factor $\omega=(2-\alpha)/\alpha$ shall give an optimal temporal convergence $\mathcal{O}(N^{-(2-\alpha)})$ with $N$ denoting the number of time mesh. Finally, ample numerical examples are delivered to illustrate our theoretical analysis and the efficiency of the suggested scheme.
\end{abstract}
\end{frontmatter}
%%
%% Start line numbering here if you want
%%
% \linenumbers
%% main text
\section{Introduction}%\label{1}
	With the development of science and technology, it has been found that normal diffusion models can not adequately describe many natural phenomena in complex dynamic systems that exhibit power-law decaying behavior \cite{Klj2005AD,Buc2015ND,Duos2019AC,Duq2019NM,Dew2020MA}. Such diffusion usually deviates from the assumption of Brownian motion so that it is called anomalous diffusion. In addition, from a probabilistic point of view, anomalous diffusion commonly follows a relation that the second moment or mean squared displacement of a particle is a nonlinear function of time $t$, i.e., $\langle\left|x(t)\right|^{2}\rangle\sim t^{\beta}$ with $\beta=[0,1)\cup (1,2]$, $x(t)$ is the trajectory of a particle or stochastic process \cite{Sor2008OT,Sor2003PT,Dew2018BP}. Anomalous diffusion phenomena are ubiquitous in natural world, such as diffusive transport of solutes in heterogeneous porous media \cite{Mer2000TR,Mem2011SM}, RNA movement in bacterial cytoplasm \cite{Goi2006PN}, animals' food-seeking \cite{Eda2007RL,Hun2010EC}, contaminants in groundwater \cite{Kij2000FS}, material with thermal memory \cite{Wol1994OI} and so on.
	
	Of the many possible models of anomalous diffusion, we shall be interested in so-called fractional diffusion, which is a nonlocal process and can be derived to model diffusive transport phenomena of solutes in heterogeneous porous media \cite{Mer2000TR,Mem2011SM}. Specificly, we concentrate on the following nonlocal diffusion model
	\begin{align}\label{equ:1.1}
		\begin{cases}
			\partial_{t}^{\alpha}u(\bm{x},t)+\kappa_{s}(-\Delta +\gamma \mathbb{I})^{s} u(\bm{x},t)=f(\bm{x},t),\quad \bm{x}\in \Omega,~~0<t\leq T, \\
			u(\bm{x}, 0)=u_{0}(\bm{x}),\quad \bm{x} \in \Omega, \\
			u(\bm{x}, t)=0,\quad \bm{x}\in \partial\Omega,~~0<t\leq T,
		\end{cases}
	\end{align}
    where $0<s<1$, $\gamma\geq 0$, $\kappa_{s}>0$ is diffusion coefficient, $\mathbb{I}$ denotes identity operator, $\Omega\subset \mathbb{R}^{d}$ ($d=1,2,3$) is a bounded open domain, and initial value $u_{0}$ and source term $f$ are known. Through the framework of a continuous time random walk under the assumption that the mean waiting time has a power-law decaying tail \cite{Mer2000TR,Mem2011SM}, the Caputo derivative of $u(\bm{x},t)$ can be derived as follows: 
	\begin{align}\label{equ:1.2}
		\partial_{t}^{\alpha}u(\bm{x},t)=\frac{1}{\Gamma(1-\alpha)} \int_{0}^{t} \frac{\partial u(\bm{x}, \tau)}{\partial \tau} \frac{d \tau}{(t-\tau)^{\alpha}},
	\end{align} 
	where $0 <\alpha<1$, and $\Gamma(\cdot)$ corresponds to Euler's Gamma function.
	
	For the fractional Laplacian $(-\Delta+\gamma \mathbb{I})^{s}$, a definition of tempered fractional Laplacian was reported in \cite{Dew2018BP} on bounded domain. In this work, we primarily consider the following spectral decomposition form:
	\begin{align} \label{equ:1.3}
		(-\Delta+\gamma \mathbb{I})^{s} u(\bm{x},t)=\sum_{j=1}^{\infty}(u,\varphi_{j})(\lambda_{j}+\gamma)^{s}\varphi_{j},
	\end{align}
	where $\left\{(\lambda_{j},\varphi_{j})\right\}_{j=1}^{\infty}$ denote the eigenpairs of classical Laplacian operator $(-\Delta)$ with homogeneous Dirichlet boundary conditions on $\Omega$. When taking $\gamma=0$ in (\ref{equ:1.3}), the fractional Laplacian can represent the infinitesimal generator of a subordinate killed Brownian motion, namely, the process that first kills Brownian motion in $\Omega$ and then subordinates it through an $s$-stable subordinator \cite{Sor2008OT,Sor2003PT,Lix2020NA,Duos2019AC}.
	
	Time-space fractional diffusion problems involving spectral fractional Laplacian have attracted considerable attention and been widely studied in pasted decade. But all along, nonlocal nature of time-space fractional differential operators brings scholars great challenges in aspects of analysis and simulation. In \cite{NoR2016AP}, a Caffarelli-Silvestre extension turned problem (\ref{equ:1.1}) with $\gamma=0$ into a quasistationary elliptic problem with dynamic boundary condition, which brought one extra dimension and a degenerate weight. The finite element method (FEM) and spectral method based on the Caffarelli-Silvestre extension were introduced in \cite{NoR2015AP,Boa2018NM,Chs2020AE}. Through Dunford-Taylor integral representation, finite element discretization using sinc quadrature for spectral fractional Laplacian was proposed in \cite{Boa2015NA,Boa2018NM,Bona2017NA,Boa2017TA}. Recently, a best uniform rational approximation was found for solving spectral fractional elliptic problems \cite{Has2020AO,Has2019TH,Has2019NS}. But it is difficult to derive a similar finite element error estimate developed in \cite[Theorem 4.1]{Has2020AO} for time-space fractional diffusion problems. In addition, there were also a family of wildly used numerical methods combining with matrix transfer technique, such as finite difference methods \cite{Ilm2005NA,Yaq2011NN,Seb2017FD,Dus2018AF}, FEMs \cite{Yaq2011NN,Buk2012AE,Mag2020FE,Szb2021FE}, and also spectral and spectral element methods \cite{Zhm2022MT}. However, the use of matrix transfer technique will encounter a thorny problem that how to treat fractional power of matrix obtained from spatial discretization, which usually leads to high computational cost and thus is rather time-consuming.
	
	To this end, it is necessary to establish concise and implementable solutions for time-space nonlocal model problems. In this paper, we shall be devoted to an efficient implementation of matrix transfer technique. To the best of our current knowledge, some current existing works employing matrix transfer technique \cite{Yaq2011NN,Ilm2005NA,Seb2017FD,Mag2020FE,Szb2021FE,Zhm2022MT} mainly proposed the numerical approximation to fractional power of matrix, and then iteration methods (e.g., Krylov subspace method) for system of equations are necessary. In \cite{Dus2018AF}, authors reported a fast second-order central difference scheme and matrix transfer technique for time-space fractional diffusion problem, which can exactly compute the solution of diffusion problem but without theoretical analysis. In addition, we also notice that there is still no literature on finite element analysis for time-space fractional diffusion problem (\ref{equ:1.1}) using fast time-stepping L1 scheme on graded temporal mesh.
	
	This paper has two aims. The first goal is to propose an efficient and accurate spatial semi-discrete scheme using linear FEM or fourth-order compact difference method (CDM) and matrix transfer technique for solving problem (\ref{equ:1.1}). Since eigen-decomposition of symmetric tri-diagonal Toeplitz matrix coming from finite element or compact difference discretization can be carried out by $1$-dimensional discrete sine transform, we extend the technique to deal with spectral fractional Laplacian on $d$-dimensional Cartesian meshes with the aids of the properties of Kronecker product and generalized multinomial theorem. The method can exactly and fastly evaluate fractional power of matrix and perform matrix-vector multiplication using $d$-dimensional discrete sine transform at arbitrary time level, and thus enjoys low computational cost and high numerical accuracy. This is significantly different from some existing reports \cite{NoR2015AP,Boa2018NM,Ilm2005NA,Yaq2011NN,Seb2017FD,Buk2012AE,Mag2020FE,Szb2021FE}. The second is to study stability and convergence analyses of full discrete scheme using fast time-stepping L1 scheme on graded temporal mesh. The provided analysis shall illustrate that a choice of graded mesh factor $\omega=(2-\alpha)/\alpha$ can yield an optimal temporal convergence $\mathcal{O}(N^{-(2-\alpha)})$ in time, where $N$ denotes the number of time intervals.
	
	The main contributions of this paper include the following: (1) efficient finite element or compact difference semi-discretization using matrix transfer technique for time-space fractional diffusion problem is established; (2) stability and error bounds of full discrete scheme based on fast time-stepping L1 method on graded time meshes are discussed; (3) fast and accurate solver for system of equations is designed by discrete sine transform at per time level; (4) commendably concise and implementation-friendly double fast algorithm is proposed.
	
	The remainder of this paper is organized as follows. In Section 2, matrix transfer technique is introduced, and spatial semi-discrete scheme of time-space fractional diffusion problem is derived. In Section 3, full discrete scheme based on a fast time-stepping L1 scheme is developed, and the corresponding stability and error analyses are addressed. Section 4 performs some numerical experiments. A brief summary is made in Section 5.
\section{Spatial discretization}\label{sec.3}
	Since there are various approximation techniques for spectral fractional Laplacian problems, we shall be committed to study the matrix transform technique \cite{Ilm2005NA} in this section.
	
	To this end, we first give a concise introduction of matrix transform technique by using linear FEM or fourth-order CDM to discretize fractional Possion problem. Then we extend the numerical algorithm to time-space fractional diffusion problem (\ref{equ:1.1}), and also derive the corresponding spatial semi-discrete scheme.
\subsection{Matrix transform technique}%Tensorial numerical methods
	Throughout this paper, let $\Omega=\prod\limits_{k=1}^{d}(a_{k},b_{k})$ be a bounded and open domain with $a_{k}<b_{k}$, denote by $N_{k}$ the number of meshgrids, take the spatial mesh size $h_{k}=\frac{b_{k}-a_{k}}{N_{k}}$, and define a finite dimensional space
	\begin{align}\notag
		V_{h_{k}}=span\left\{\phi_{0}^{(k)},\phi_{1}^{(k)},...,\phi_{N_{k}}^{(k)}\right\},\quad~ k = 1,2,...,d,
	\end{align}
	where $\phi_{j}^{(k)}$ corresponds to piecewise linear basis function for $j=0,1,...,N_{k}$. Then the $d$-dimensional tensorial finite element space $V_{h}^{d}$ can be given by 
	\begin{align}\notag
	V_{h}^{d}=V_{h_{1}}\otimes V_{h_{2}}\otimes \cdots \otimes V_{h_{d}},\quad d=1,~2,~3,
	\end{align}
	where $\otimes$ represents Kronecker product. In addition, we also define an index set
	\begin{align}\notag
		\mathcal{T}_{d}=\left\{i=(i_{1},i_{2},...,i_{d})~|~1\leq i_{k}\leq N_{k}-1,~k=1,2,...,d\right\}.
	\end{align}
	
	We now consider tensorial linear FEM or fourth-order CDM for fractional Possion problem
	\begin{align}\label{equ:2.1}
		\begin{cases}
			(-\Delta+\gamma\mathbb{I})^{s}u(\bm{x})=f(\bm{x}),\quad \bm{x}\in \Omega, \\
			u(\bm{x})=0,\quad \bm{x}\in \partial\Omega.
		\end{cases}
	\end{align}
	By matrix transform technique, (\ref{equ:2.1}) is formulated as a matrix-vector multiplication form:
	\begin{align}\label{equ:2.2}
		(M^{-1}S+\gamma I)^{s}\bm{u}=K\bm{b},
	\end{align}
	where $I$, $S$, $M$ and $\bm{b}$ denote identity matrix, stiffness matrix, mass matrix and right hand side vector, respectively. Here, we let
	\begin{align}\notag
		K=K_{1}\otimes\cdots\otimes K_{k-1}\otimes K_{k}\otimes K_{k+1}\otimes\cdots \otimes K_{d},
	\end{align}
	with $K_{k}=I_{k}$ or $M_{k}^{-1}$ using finite element discretization, while $K_{k}=I_{k}$ by CDM. The explicit $I_{k}$ and $M_{k}$ will be defined later. Most notably, we also have
	\begin{align}\notag
		S=\sum_{k=1}^{d}\left(M_{1}\otimes\cdots\otimes M_{k-1}\otimes A_{k}\otimes M_{k+1}\otimes\cdots \otimes M_{d}\right),
	\end{align}
	as well as
	\begin{align}\notag
		M=M_{1}\otimes\cdots\otimes M_{k-1}\otimes M_{k}\otimes M_{k+1}\otimes\cdots \otimes M_{d},\quad I=I_{1}\otimes\cdots\otimes I_{k-1}\otimes I_{k}\otimes I_{k+1}\otimes\cdots \otimes I_{d},
	\end{align}
	where $I_{k}$ is $(N_{k}-1)$ by $(N_{k}-1)$ identity matrix,
	\begin{align}\notag
		A_{k}=\left(\begin{array}{llllll}
			a     &   b   & 0    &   0   & \cdots & 0\\
			b     &   a   & b     &   0   & \cdots & 0\\
			\vdots     & \ddots & \ddots &  \ddots & \ddots &\vdots\\
			0 & \cdots &   0  & b      & a    & b\\
			0 & \cdots &   0  & 0     & b     & a
		\end{array}\right),\quad and~~
		M_{k}=\left(\begin{array}{llllll}
			c     &     d  &    0   &   0    & \cdots &   0\\
			d      &    c  &    d   &   0    & \cdots &   0\\
			\vdots     & \ddots & \ddots & \ddots & \ddots &  \vdots\\
			0     & \cdots &    0   &   d    &    c   &   d\\
			0     & \cdots &    0   &   0    &    d   &   c
		\end{array}\right),
	\end{align}
    here, $a=\frac{2}{h_{k}}$, $b=\frac{-1}{h_{k}}$, $c=\frac{4h_{k}}{6}$ and $d=\frac{h_{k}}{6}$ by FEM, and yet $a=\frac{2}{h_{k}^{2}}$, $b=\frac{-1}{h_{k}^{2}}$, $c=\frac{10}{12}$ and $d=\frac{1}{12}$ by CDM. It is noteworthy that $(N_{k}-1)$ by $(N_{k}-1)$ matrices $A_{k}$ and $M_{k}$ are symmetric and tri-diagonal, their eigenvalues can be expressed as
	\begin{align}\notag
		\lambda_{i}^{\left(s_{k}\right)}=a-2\left|b\right|\cos\left(\frac{i\pi}{N_{k}}\right),\quad and~~~
		\lambda_{i}^{\left(m_{k}\right)}=c+2\left|d\right|\cos\left(\frac{i\pi}{N_{k}}\right),
	\end{align}
	respectively, and the eigenvectors of both of them are $\left(P_{k}\right)_{ij}=\sqrt{\frac{2}{N_{k}}}\sin\left(\frac{ij\pi}{N_{k}}\right)$ for $i,j=1,2,\dots,N_{k}-1$. For more details, see \cite{Elj1953TC,Ler1955FD,Grr1969AC}. Therefore, the eigen-decompositions of $A_{k}$ and $M_{k}$ are separately as follows:
	\begin{align}\label{equ:2.3}
		A_{k}=P_{k}\Lambda_{s_{k}}P_{k}^{-1},\quad and~~
		M_{k}=P_{k}\Lambda_{m_{k}}P_{k}^{-1},
	\end{align}
	where the diagonal matries
	\begin{align}\notag
	\Lambda_{s_{k}}=diag\left(\lambda_{1}^{(s_{k})},\lambda_{2}^{(s_{k})},\dots,\lambda_{N_{k}-1}^{(s_{k})}\right),\quad and~~~
	\Lambda_{m_{k}}=diag\left(\lambda_{1}^{(m_{k})},\lambda_{2}^{(m_{k})},\dots,\lambda_{N_{k}-1}^{(m_{k})}\right).
	\end{align}
	Also, it follows that $P_{k}=P_{k}^{T}=P_{k}^{-1}$ for $k=1,2,...,d$.
	
	Below, we shall extend the above numerical algorithm to time-space fractional diffusion problem (\ref{equ:1.1}), which will be stated in the next subsection.	
\subsection{Spatial semi-discretization}
	In this part, we primarily present a commendably concise and implementation-friendly spatial numerical algorithm for discreting fractional diffusion problem (\ref{equ:1.1}).
	
	Based on the discretization of fractional Possion problem (\ref{equ:2.1}), we immediately have the following matrix-vector multiplication form of (\ref{equ:1.1}):
	\begin{align}\label{equ:2.4}
		\left(\partial_{t}^{\alpha}+\kappa_{s}(M^{-1}S+\gamma I)^{s}\right)\bm{u}(t)=K\bm{b}(t),
	\end{align}
	where $M$, $S$ and $K$ are given by (\ref{equ:2.2}), and $\bm{b}$ denotes a time-dependent right hand side vector function. The next work aims to deal with the fractional power $(M^{-1}S+\gamma I)^{s}$ appearing in (\ref{equ:2.4}). 
	
	In terms of the inverse properties of Kronecker products and its matrix multiplication, we have
	\begin{align}\notag
		M^{-1}S&=\sum_{k=1}^{d}\left((M_{1}^{-1}M_{1})\otimes\cdots\otimes (M_{k-1}^{-1}M_{k-1})\otimes(M_{k}^{-1}A_{k})\otimes (M_{k+1}^{-1}M_{k+1})\otimes\cdots \otimes(M_{d}^{-1}M_{d})\right)
		\notag\\
		&=\sum_{k=1}^{d}\left(I_{1}\otimes\cdots\otimes I_{k-1}\otimes(M_{k}^{-1}A_{k})\otimes I_{k+1}\otimes\cdots \otimes I_{d}\right)
		\notag.
	\end{align}
	Therefore, we can obtain
	\begin{align}\label{equ:2.5}
		(M^{-1}S+\gamma I)&=
		\sum_{k=1}^{d}\left(I_{1}\otimes\cdots\otimes I_{k-1}\otimes(M_{k}^{-1}A_{k})\otimes I_{k+1}\otimes\cdots \otimes I_{d}\right. \notag \\
		&+\left.\frac{\gamma}{d}I_{1}\otimes\cdots\otimes I_{k-1}\otimes I_{k}\otimes I_{k+1}\otimes\cdots \otimes I_{d}\right) \notag\\
		&=\sum_{k=1}^{d}\left(I_{1}\otimes\cdots\otimes I_{k-1}\otimes(M_{k}^{-1}A_{k}+\frac{\gamma}{d}I_{k})\otimes I_{k+1}\otimes\cdots \otimes I_{d}\right).
	\end{align}
	Moreover, leveraging eigen-decomposition (\ref{equ:2.3}), we derive
	\begin{align}\label{equ:2.6}
		T_{k}:=(M_{k}^{-1}A_{k}+\frac{\gamma}{d}I_{k})
		=P_{k}\Lambda_{k}P_{k}^{-1},
	\end{align}
	where $\Lambda_{k}:=diag\left(\lambda_{1}^{(k)},\lambda_{2}^{(k)},\dots,\lambda_{N_{k}-1}^{(k)}\right)$ with $\lambda_{i}^{(k)}=\left(\frac{\lambda_{i}^{(s_{k})}}{\lambda_{i}^{(m_{k})}}+\frac{\gamma}{d}\right)$ for $i=1,2,\dots,N_{k}-1$.
	
	Combining (\ref{equ:2.5}) and (\ref{equ:2.6}), and using generalized multinomial theorem (see \cite{Gaj2010DM,No2022TR,AliGM2022}) yield
	\begin{align}\label{equ:2.7}
		\left(M^{-1}S+\gamma I\right)^{s}&=\left(\sum_{k=1}^{d}\left(I_{1}\otimes\cdots\otimes I_{k-1}\otimes T_{k}\otimes I_{k+1}\otimes\cdots \otimes I_{d}\right)
		\right)^{s} \notag\\
		&=\sum_{n_{d-1}=0}^{\infty}\sum_{n_{d-2}=0}^{n_{d-1}}\cdots\sum_{n_{1}=0}^{n_{2}}\binom{s}{n_{d-1}}\binom{n_{d-1}}{n_{d-2}}\cdots\binom{n_{2}}{n_{1}}
		\left(I_{1}\otimes\cdots\otimes I_{d-1}\otimes T_{d}\right)^{s-n_{d-1}}\notag\\
		&\times\left(I_{1}\otimes\cdots I_{d-2}\otimes T_{d-1} \otimes I_{d}\right)^{n_{d-1}-n_{d-2}}\cdots
		\left(T_{1}\otimes I_{2}\otimes\cdots\otimes I_{d}\right)^{n_{1}}\notag\\
		&=\sum_{n_{d-1}=0}^{\infty}\sum_{n_{d-2}=0}^{\infty}\cdots\sum_{n_{1}=0}^{\infty}\binom{s}{n_{1},n_{2},\ldots,n_{d-1}}\left(T_{1}\otimes I_{2}\otimes\cdots\otimes I_{d}\right)^{n_{1}}\cdots\notag\\ 
		&\times
		\left(I_{1}\otimes\cdots I_{d-2}\otimes T_{d-1} \otimes I_{d}\right)^{n_{d-1}}\left(I_{1}\otimes\cdots\otimes I_{d-1}\otimes T_{d} \right)^{s-n_{1}-\cdots-n_{d-1}},
	\end{align}
	where the generalized multinomial coefficient
	\begin{align}\notag
	\binom{s}{n_{1},n_{2},\ldots,n_{d-1}}:=\binom{s}{n_{1}+n_{2}+\cdots+n_{d-1}}\binom{n_{1}+n_{2}+\cdots+n_{d-1}}{n_{1}+n_{2}+\cdots+n_{d-2}}\cdots\binom{n_{1}+n_{2}}{n_{1}}.
	\end{align}
    %To provide transparency and acknowledge the origins of these techniques, we recommend referring to \cite{Gaj2010DM,No2022TR,AliGM2022}.
	
	In what follows, the treatment means of (\ref{equ:2.7}) will be coincidence with \cite{Dus2018AF}. So we directly have 
	\begin{align}\label{equ:2.8}
		\left(M^{-1}S+\gamma I\right)^{s}\bm{u}(t)=P_{1}\textcircled{x}_{1}\cdots\textcircled{x}_{d-1}P_{d}\textcircled{x}_{d}\left(H_{d}\odot
		\left(P_{1}^{-1}\textcircled{x}_{1}\cdots\textcircled{x}_{d-1}P_{d}^{-1}\textcircled{x}_{d}U(t)\right)
		\right),
	\end{align}
	where $H_{d}$ is a $d$-dimensional array denoted by
	\begin{align}\label{equ:2.9}
		\left(H_{d}\right)_{i_{1},i_{2},\ldots,i_{d}}=\left(\sum_{k=1}^{d}\lambda_{i_{k}}^{(k)}\right)^{s}, \quad\left(i_{1}, i_{2}, \ldots, i_{d}\right) \in \mathcal{T}_{d}.
	\end{align}
	Further, $\odot$ represents Hadamard array product, $U(t)=\left\{\bm{u}_{i}(t)\right\}_{(N_{1}-1)\times\cdots\times(N_{d}-1)}$ is a $d$-dimensional array for $i \in \mathcal{T}_{d}$, and the notation $\textcircled{x}_{k}$ is defined by
	\begin{align}\notag
	 	\left(W\textcircled{x}_{k}R\right)_{i_{1},i_{2},\ldots,i_{d}}:=\sum_{j=1}^{N_{k}-1}W_{i_{k},j}R_{i_{1},\ldots,i_{k-1},j,i_{k+1}\ldots,i_{d}},
	\end{align}
	where $W$ corresponds to an $\left(N_{k}-1\right)$ by $\left(N_{k}-1\right)$ matrix, and $R$ denotes a $d$-dimensional array.
	
	Next, we give an explicit computation of right hand side vector. Analogous to (\ref{equ:2.8}), we shall have
	\begin{align}\label{equ:2.10}
		K\bm{b}(t)&=\left(K_{1}\otimes\cdots\otimes K_{k-1}\otimes K_{k}\otimes K_{k+1}\otimes\cdots \otimes K_{d}\right)\bm{b}(t)\notag\\
		&=P_{1}\textcircled{x}_{1}\cdots\textcircled{x}_{d-1}P_{d}\textcircled{x}_{d}\left(V_{d}\odot
		\left(P_{1}^{-1}\textcircled{x}_{1}\cdots\textcircled{x}_{d-1}P_{d}^{-1}\textcircled{x}_{d}F(t)\right)
		\right),
	\end{align}
	where $V_{d}$ is a $d$-dimensional array with entry
	\begin{align}\label{equ:2.11}
		\left(V_{d}\right)_{i_{1},i_{2},\ldots,i_{d}}=\prod_{k=1}^{d}\bar{\lambda}_{i_{k}}^{(k)}, \quad\left(i_{1}, i_{2}, \ldots, i_{d}\right) \in \mathcal{T}_{d},
	\end{align}
	where $\bar{\lambda}_{i_{k}}^{(k)}$ means the $i_{k}$'th eigenvalue of $K_{k}$, and $F(t)=\left\{\bm{b}_{i_{1},i_{2},\ldots,i_{d}}(t)\right\}_{(N_{1}-1)\times\cdots\times(N_{d}-1)}$ is a $d$-dimensional array given by
	\begin{align}\label{equ:2.12}
		\bm{b}_{i_{1},i_{2},\ldots,i_{d}}(t)=f(x_{1_{i_{1}}},x_{2_{i_{2}}},\ldots,x_{d_{i_{d}}},t),\quad or~~~ \bm{b}_{i_{1},i_{2},\ldots,i_{d}}(t)=\left(f(\cdot,t),\prod\limits_{k=1}^{d}\phi_{i_{k}}^{(k)}(\cdot)\right),
	\end{align}
	where $\left(\cdot,\cdot\right)$ denotes $L^{2}$ inner product which will be defined later. It should be pointed out by (\ref{equ:2.12}) that the former is obtained using finite element discretization with assumption $f(t)\approx f_{h}(t)\in V_{h}^{d}$ or employing CDM directly, at this time one shall choose $K=I$; however, the latter will appear in finite element case, then we have to take $K=M^{-1}$.
	
	In addition, it follows from \cite{Jul2015FE,Zhl2016FH,Dus2018AF} that
	\begin{align}\label{equ:2.13}
		P_{1}\textcircled{x}_{1}\cdots\textcircled{x}_{d-1}P_{d}\textcircled{x}_{d}R=\mathscr{D}_{d}\left(R\right),\quad and ~~~
		P_{1}^{-1}\textcircled{x}_{1}\cdots\textcircled{x}_{d-1}P_{d}^{-1}\textcircled{x}_{d}R=\mathscr{D}^{-1}_{d}\left(R\right),
	\end{align}
	where $\mathscr{D}_{d}$ and $\mathscr{D}_{d}^{-1}$ represent $d$-dimensional discrete sine transform and its inverse transform, respectively, and $R$ is a $d$-dimensional array. 
	
	Hence, using (\ref{equ:2.4}), (\ref{equ:2.8}), (\ref{equ:2.10}) and (\ref{equ:2.13}), the semi-discrete scheme of problem (\ref{equ:1.1}) can be written as
	\begin{align}\label{equ:2.14}
		\partial_{t}^{\alpha}U(t)+\mathscr{D}_{d}\left(\kappa_{s}H_{d}\odot \mathscr{D}_{d}^{-1}\left(U(t)\right)\right)=\mathscr{D}_{d}\left(V_{d}\odot \mathscr{D}_{d}^{-1}\left(F(t)\right)\right),\quad t\in(0,T].
	\end{align}
	
	\begin{remark}\label{re:2.1} For the fractional Possion problem (\ref{equ:2.1}), its numerical solution becomes
		\begin{align}\notag
			U=\mathscr{D} _{d}\left(L_{d} \odot \left(V_{d}\odot\mathscr{D}_{d}^{-1}(F)\right)\right),
		\end{align}
		where $F=\left\{\bm{b}_{i}\right\}_{(N_{1}-1)\times\cdots\times(N_{d}-1)}$ is a $d$-dimensional array for $i \in \mathcal{T}_{d}$, and $L_{d}$ is defined by
		\begin{align}\notag
			\left(L_{d}\right)_{i_{1},i_{2},\ldots, i_{d}}=1 /\left(\kappa_{s}H_{d}\right)_{i_{1},i_{2}, \ldots, i_{d}}, \quad\left(i_{1}, i_{2}, \ldots, i_{d}\right) \in \mathcal{T}_{d}.
		\end{align}
		In particular, if taking $f(\bm{x})\approx f_{h}(\bm{x})\in V_{h}^{d}$, then
		\begin{align}\notag
			U=\mathscr{D} _{d}\left(L_{d} \odot \mathscr{D}_{d}^{-1}(F)\right),
		\end{align}
		where $F=\left\{f(x_{1_{i_{1}}},x_{2_{i_{2}}},\ldots,x_{d_{i_{d}}})\right\}_{(N_{1}-1)\times\cdots\times(N_{d}-1)}$ denotes a $d$-dimensional array.
	\end{remark}
	
	 From Remark \ref{re:2.1}, it should be mentioned that the fast algorithm can reduce the computational cost from $\mathcal{O}\left(\mathcal{M}^{d+1}\right)$ by a direct matrix-vector multiplication to $\mathcal{O}\left(\mathcal{M}^{d}\log_{2}\mathcal{M}\right)$ with $\mathcal{M}=\max\limits_{1 \leq k \leq d}{N_{k}}$.
	
	\begin{remark}\label{re:2.2} If $M=K=I$ in the fourth-order CDM, then (\ref{equ:2.4}) or (\ref{equ:2.14}) reduces to the second-order central finite difference semi-discrete scheme \cite{Dus2018AF} for (\ref{equ:1.1}) in space.
	\end{remark}
	
	For the present algorithm above, due to the uses of discrete sine transform and its inverse transform, we can exactly evaluate fractional power of matrix (\ref{equ:2.7}) and perform the matrix-vector multiplication in (\ref{equ:2.14}), which is different from some existing techniques \cite{NoR2015AP,Boa2018NM,Ilm2005NA,Yaq2011NN,Seb2017FD,Buk2012AE,Mag2020FE,Szb2021FE} with high computational cost and low numerical accuracy.
	
	 Hereto, we have completed the introduction on the spatial numerical algorithm, which is both concise and implementable. In the next section, we will address the implementation and convergence analysis of full discrete scheme which is based on a fast time-stepping L1 method \cite{Jis2017FE,Yay2017FE,Huy2022AN,Huy2022EE} in time. 
\section{Implementation and analysis of full discrete scheme}
	Though a fast Fourier algorithm has been presented in space, it is still difficult for us to numerically simulate the high dimensional time-space fractional diffusion problem (\ref{equ:1.1}) in temporal direction. Hence, it is critical to choose a discrete technique with low memory and computational cost for problem (\ref{equ:1.1}) involving time fractional Caputo derivative (\ref{equ:1.2}).
	
	In this part, we first introduce a fast time-stepping L1 scheme for time discretization, and induce the full discrete scheme of (\ref{equ:1.1}). Further, we shall give convergence analysis of full discrete scheme.
\subsection{Fast time-stepping method}
	Let $N$ be the number of temporal meshes, and $\mathcal{I}^{1}$ the piecewise linear Lagrange interpolation operator. Given temporal graded mesh $t_{n}=T\left(n/N\right)^{\omega}$ with $\omega\geq 1$ and $n=1,2,...,N$, we take $t=t_{n}$ such that 
	\begin{align}\notag
		\partial_{t}^{\alpha} \bm{u}\left(t_{n}\right)=\frac{1}{\Gamma(1-\alpha)} \int_{0}^{t_{n}} \frac{d\bm{u}(\zeta)}{d \zeta}\left(t_{n}-\zeta\right)^{-\alpha}d\zeta,\quad n=1,2,...,N. \notag
	\end{align}% on interval $[t_{n-1},t_{n}]$
	
	Interpolating linearly $\bm{u}(\zeta)$ with respect to variable $\zeta$ on $[t_{k-1},t_{k}]$ for $k=1,2,...,n$, and using the integration by parts on $[0,t_{n-1}]$, we have
	\begin{align}\label{equ:3.1}
		\partial_{t}^{\alpha} \bm{u}\left(t_{n}\right) &\approx\frac{1}{\Gamma(1-\alpha)} \int_{t_{n-1}}^{t_{n}} \frac{d\mathcal{I}^{1}\bm{u}(\zeta)}{d \zeta}\left(t_{n}-\zeta\right)^{-\alpha}d\zeta
		+\frac{1}{\Gamma(1-\alpha)} \int_{0}^{t_{n-1}} \frac{d\mathcal{I}^{1}\bm{u}(\zeta)}{d\zeta}\left(t_{n}-\zeta\right)^{-\alpha}d\zeta \notag \\
		&=\frac{\bm{u}(t_{n})-\alpha \bm{u}(t_{n-1})-(1-\alpha)\Delta t_{n}^{\alpha}t_{n}^{-\alpha}\bm{u}(t_{0})}{\tau_{n}}
		+\frac{1}{\Gamma(-\alpha)}\int_{0}^{t_{n-1}} \frac{\mathcal{I}^{1}\bm{u}(\zeta)}{\left(t_{n}-\zeta\right)^{1+\alpha}}d\zeta\notag\\
		&=:L^{\alpha}\bm{u}(t_{n})+H^{\alpha}\bm{u}(t_{n}),
	\end{align}
	where $\tau_{n}=\Delta t_{n}^{\alpha}\Gamma(2-\alpha)$, $L^{\alpha}\bm{u}(t_{n})$ and $H^{\alpha}\bm{u}(t_{n})$ are called local part and history part, respectively.

	In order to develop fast time-stepping L1 method, we now introduce a sum of exponentials \cite{Jis2017FE,Yay2017FE,Huy2022AN,Huy2022EE} to approximate the kernel $\left(t_{n}-\zeta\right)^{-1-\alpha}$. In this work, we primarily consider the approximation which is based on the trapezoidal rule on a real line. Namely, using the definition of Euler integral representation of Gamma function, a change of variable $s=e^{y}$ and $\Gamma(1-\alpha)\Gamma(\alpha)=\frac{\pi}{\sin(\pi\alpha)}$ gives
	\begin{align}\label{equ:3.2}
		\frac{t^{-\alpha-1}}{\Gamma(-\alpha)}&=\frac{1}{\Gamma(-\alpha)\Gamma(\alpha+1)}\int_{0}^{\infty}s^{\alpha}e^{-st}ds=\int_{-\infty}^{\infty}g(y,t,\alpha)dy,
	\end{align}
	where $g(y,t,\alpha)=-\frac{\sin(\pi\alpha)}{\pi}e^{(\alpha+1)y}e^{-e^{y}t}$ which decays exponentially as $|y|\to \infty$. Therefore, we can use the exponentially convergent trapezoidal rule \cite{Trl2014TE} to approximate the integral
	\begin{align}\label{equ:3.3}
		\int_{-\infty}^{\infty}g(y,t,\alpha)dy\approx\Delta y\sum_{j=-\infty}^{\infty}g(j\Delta y,t,\alpha)=\sum_{j=1}^{Q} w_{j} e^{\xi_{j} t}+O\left(\epsilon t^{-\alpha-1}\right), \quad t \in\left[t_{1}, T\right],
	\end{align}
	where $Q>1$ is a positive integer, $w_j=-\frac{\sin (\alpha \pi)}{\pi} \Delta y e^{(1+\alpha) y_j}$, $\xi_j=-e^{y_j}$, $y_j=y_{\min }+(j-1)\Delta y$, $\Delta y=\frac{y_{\max }-y_{\min }}{Q-1}$, and $\epsilon>0$ is a given precision satisfying
	\begin{align}\label{equ:3.4}
		\left|\int_{-\infty}^{a}g(y,t,\alpha)dy\right|+\left|\int_{b}^{\infty}g(y,t,\alpha)dy\right|\leq 2\epsilon t^{-1-\alpha},\quad a<0<b.
	\end{align}
	The use of (\ref{equ:3.4}) determines $y_{\min}=(1+\alpha)^{-1}\ln\left(\epsilon\right)-\ln(T)\geq a$ and $y_{\max}=\ln\left(\frac{-\ln \left(\epsilon\right)+(1+\alpha) \ln \left(t_{1}\right)}{0.5t_{1}}\right)\leq b$. For some detailed derivations, we refer reader to \cite[Appendix C]{Gul2019EM}.
	
	Indeed, the above inequality (\ref{equ:3.4}) is obvious. A change of variable leads to
	\begin{align}\notag
		\left|\int_{b}^{\infty}g(y,t,\alpha)dy\right|\leq t^{-\alpha-1}\Gamma(\alpha+1,te^{b}),\quad and~~ \left|\int_{-\infty}^{a}g(y,t,\alpha)dy\right|\leq t^{-\alpha-1}\gamma(\alpha+1,te^{a}),
	\end{align}
	where $\Gamma(\cdot,\cdot)$ and $\gamma(\cdot,\cdot)$ denote upper and lower incomplete Gamma functions, respectively. For larger $b$ ($b>0$) and smaller $a$ ($a<0$), $\Gamma(\alpha+1,te^{b})$ and $\gamma(\alpha+1,te^{a})$ decay exponentially to 0 so that it is feasible to take $\epsilon$ to be as follows:
	\begin{align}\label{equ:3.5}
		\epsilon=\max\left\{\Gamma(\alpha+1,te^{b}),\gamma(\alpha+1,te^{a})\right\}.
	\end{align}
	Thus we obtain (\ref{equ:3.4}).
	
	Hence, from (\ref{equ:3.2}) and (\ref{equ:3.3}), we immediately have
	\begin{align}\label{equ:3.6}
		H^{\alpha}\bm{u}(t_{n})=\sum_{j=1}^{Q} w_{j} e^{\xi_{j}\Delta t_{n}} Y_{j}\left(t_{n-1}\right),%+O(\epsilon\Delta t)
	\end{align}
	where $Y_{j}(t)=\int_{0}^{t} e^{\xi_{j}(t-\zeta)}\mathcal{I}^{1}\bm{u}(\zeta)d\zeta$ satisfies 
	\begin{align}\notag
	(Y_{j})^{\prime}(t)=\xi_{j} Y_{j}(t)+\mathcal{I}^{1}\bm{u}(t) \quad with~~ Y_{j}(0)=0, 
	\end{align}
	which can be exactly solved by 
	\begin{align}\notag
	Y_{j}\left(t_{i}\right)=\kappa_{1} Y_{j}\left(t_{i-1}\right)+\kappa_{2} \bm{u}(t_{i-1})+\kappa_{3} \bm{u}(t_{i})\quad with~~ Y_{j}\left(t_{0}\right)=0,\quad i\geq1,
	\end{align}
	where $\kappa_{1}=e^{z}$, $\kappa_{2}=\frac{1}{\xi_{j}}\left(e^{z}-1-\frac{e^{z}-z-1}{z}\right)$, $\kappa_{3}=\frac{1}{\xi_{j}} \frac{e^{z}-z-1}{z}$ and $z=\xi_{j}\Delta t_{i}$.
	
	Combining (\ref{equ:3.1}) and (\ref{equ:3.6}), $\partial_{t}^{\alpha} \bm{u}(t)$ at $t=t_{n}$ can be approximated by
	\begin{align}\label{equ:3.7}
	\partial_{t}^{\alpha} \bm{u}\left(t_{n}\right)=\bar{\partial}_{t}^{\alpha} \bm{u}\left(t_{n}\right)+r^{n},
	\end{align}
	where $r^{n}$ denotes truncation error, and the discrete differential operator $\bar{\partial}_{t}^{\alpha}$ is defined by
	\begin{align}\notag
		\bar{\partial}_{t}^{\alpha} \bm{u}\left(t_{n}\right):=\frac{\bm{u}(t_{n})-\alpha\bm{u}(t_{n-1})-(1-\alpha)\Delta t_{n}^{\alpha}t_{n}^{-\alpha}\bm{u}(t_{0})
			}{\tau_{n}}+\sum_{j=1}^{Q} w_{j} e^{\xi_{j}\Delta t_{n}} Y_{j}\left(t_{n-1}\right),\quad n=1,2,...,N,
	\end{align}
	where the above sum term will vanish for $n=1$.

	From (\ref{equ:2.4}) and (\ref{equ:3.7}), we shall have
	\begin{align}\label{equ:3.8}
		\left(1+\tau_{n} \kappa_{s}(M^{-1}S+\gamma I)^{s}\right)\bm{u}(t_{n})&=\bm{g}(t_{n})+\tau_{n} K\bm{b}(t_{n})-r^{n},
	\end{align}
	where $\bm{g}(t_{n})$ is defined by
	\begin{align}\notag
		\bm{g}(t_{n}):=\alpha\bm{u}(t_{n-1})+
		(1-\alpha)\Delta t_{n}^{\alpha}t_{n}^{-\alpha}\bm{u}(t_{0})-\tau_{n}\sum_{j=1}^{Q} w_{j} e^{\xi_{j}\Delta t_{n}} Y_{j}\left(t_{n-1}\right).
	\end{align}
	
	Now let $\bm{u}^{n}$ be an approximation to $\bm{u}(t_{n})$, and $\bm{g}_{i}^{n}$ the approximation to $\bm{g}_{i}(t_{n})$ for $i \in \mathcal{T}_{d}$. From (\ref{equ:2.8}), (\ref{equ:2.10}), (\ref{equ:2.13}) and (\ref{equ:3.8}), we obtain the following fully discrete scheme
	\begin{align}\label{equ:3.9}
		\mathscr{D}_{d}\left((E_{d}+\tau_{n}\kappa_{s}H_{d})\odot \mathscr{D}_{d}^{-1}\left(U^{n}\right)\right)=G_{d}^{n}+\tau_{n}\mathscr{D}_{d}\left(V_{d}\odot \mathscr{D}_{d}^{-1}\left(F^{n}\right)\right),\quad n=1,2,...,N,
	\end{align}
	where $E_{d}$ is a $d$-dimensional array with entry $(E_{d})_{i}=1$ for $i \in \mathcal{T}_{d}$, $H_{d}$ and $V_{d}$ are given in (\ref{equ:2.9}) and (\ref{equ:2.11}), respectively. In addition, $U^{n}=\left\{\bm{u}_{i}^{n}\right\}_{(N_{1}-1)\times\cdots\times(N_{d}-1)}$, $F^{n}=\left\{\bm{b}_{i}(t_{n})\right\}_{(N_{1}-1)\times\cdots\times(N_{d}-1)}$ and $G_{d}^{n}=\left\{\bm{g}_{i}^{n}\right\}_{(N_{1}-1)\times\cdots\times(N_{d}-1)}$ are three $d$-dimensional arrays for $i \in \mathcal{T}_{d}$.
	
	By taking inverse and forward discrete sine transform on both sides of (\ref{equ:3.9}) successively, we then have
	\begin{align}\label{equ:3.10}
		U^{n}=\mathscr{D}_{d}\left(L_{d}\odot\left(G_{d}^{n}+\tau_{n}(V_{d}\odot \mathscr{D}_{d}^{-1}(F^{n}))\right)\right),\quad n=1,2,...,N,
	\end{align}
	where $L_{d}$ denotes a $d$-dimensional array defined by
	\begin{align}\notag
		\left(L_{d}\right)_{i_{1},i_{2},...,i_{d}}:=\frac{1}{1+\tau_{n} k_{s}(H_{d})_{i_{1},i_{2},...,i_{d}}},  \quad\left(i_{1}, i_{2}, \ldots, i_{d}\right) \in \mathcal{T}_{d}.
	\end{align}
	
	Numerical solutions obtained by direct L1 scheme possess heavy history dependence, namely, the solutions of the current time step $t_{n}$ always depends on the solutions of all previous time steps $t_{k}$ for $k=0,1,...,n-1$. However, the fast scheme (\ref{equ:3.10}) shows that the solution of the current time step $t_{n}$ indirectly depends on the solutions of time steps $t_{n-1}$, $t_{n-2}$ and $t_{0}$. For overall computational complexity, the present double fast algorithm reduces computational cost from $\mathcal{O}\left(N^{2}\mathcal{M}^{d+1}\right)$ to $\mathcal{O}\left(NQ\mathcal{M}^{d}\log_{2}\mathcal{M}\right)$, where $Q\ll N$.
	
	\begin{remark}\label{re:3.1}
	As $\alpha \rightarrow 1^{-}$, the numerical solution in (\ref{equ:3.10}) reduces to
	\begin{align}\notag
		U^{n}=\mathscr{D}_{d}\left(L_{d}\odot\left(\mathscr{D}_{d}^{-1}(U^{n-1})+\tau_{n}(V_{d}\odot \mathscr{D}_{d}^{-1}(F^{n}))\right)\right),\quad n=1,2,...,N,
	\end{align}
	where $\tau_{n}=\Delta t_{n}$, i.e., a standard first-order backward difference scheme is used for space fractional diffusion problem.
	\end{remark}
\subsection{Convergence analysis}
	In this subsection, we are devoted to studying stability and convergence analysis of full discrete scheme (\ref{equ:3.9}) based on fast time-stepping L1 scheme. To this end, we first recall some basic function spaces, and then introduce Mittag-Leffler function and standard projection estimate. Afterwards, we also give several lemmata which play important roles in the subsequent analysis.
	
	Now, it remains to recall some basic function spaces. Let $H^{l}(\Omega)$ and $H_{0}^{l}(\Omega)$ denote Sobolev spaces with $l\geq-1$, and $\|\cdot\|_{l}$ and $\left(\cdot,\cdot\right)_{l}$ represent the corresponding $H^{l}$-norm and inner product, respectively. If $l=0$, we shall omit subscript $l$. For $s \geq-1$, $\dot{H}^{s}(\Omega) \subset H^{-1}(\Omega)$ denotes a Hilbert space equipped with norm
	\begin{align}\notag
		\|v\|_{\dot{H}^{s}(\Omega)}=\sqrt{\sum_{j=1}^{\infty}\left(\lambda_{j}+\gamma\right)^{s}\left(v,\varphi_{j}\right)^{2}},\quad \gamma\geq 0.
	\end{align}
	Given a Banach space $X$ and arbitrary $p \geq 1$, we can define
	\begin{align}\notag
		L^{p}(0, T ; X)=\left\{v(t) \in X \text { for a.e. } t \in(0, T) \text { and }\|v\|_{L^{p}(0, T ; X)}<\infty\right\},
	\end{align}
	and its norm $\|\cdot\|_{L^{p}(0, T; X)}$ is induced by
	\begin{align}\notag
		\|v\|_{L^{p}(0, T ; X)}= 
		\begin{cases}
			\left(\int_{0}^{T}\|v(t)\|_{X}^{p} d t\right)^{1 / p}, & p \in[1, \infty), \\ \operatorname{ess} \sup\limits_{t \in(0, T)}\|v(t)\|_{X}, & p=\infty .\end{cases}
	\end{align}
	
	To derive the mild solution of (\ref{equ:1.1}), it is essential to introduce Mittag-Leffler function \cite[formula 1.56]{Poi1999FD}
	\begin{align}\notag
		E_{\alpha, \beta}\left(z\right)=\sum_{k=0}^{\infty} \frac{z^{k}}{\Gamma\left(k\alpha+\beta\right)},\quad z\in \mathbb{C},
	\end{align}
	where $\beta\in \mathbb{R}$ and $0<\alpha<1$. Here, we also give the Fourier expansions of the data $u_{0}$, $f$ and $u$. Namely, let $u(t):=\sum\limits_{j=1}^{\infty}(u(\cdot,t),\varphi_{j})\varphi_{j}$, $f(t):=\sum\limits_{j=1}^{\infty}(f(\cdot,t),\varphi_{j})\varphi_{j}$ and $u_{0}:=\sum\limits_{j=1}^{\infty}(u_{0}(\cdot),\varphi_{j})\varphi_{j}$.
	
	Hence, using method of separation of variables and resorting to the Mittag-Leffler function $E_{\alpha, \beta}$, the mild solution $u$ of problem (\ref{equ:1.1}) can be represented as
	\begin{align}\label{equ:3.11}
		u(t)=E(t)u_{0}+\int_{0}^{t} \bar{E}(t-\tau) f(\tau) d\tau,
	\end{align}
	where the solution operator
	\begin{align}\notag
		E(t)u_{0}=\sum_{j=1}^{\infty}E_{\alpha, 1}\left(-\kappa_{s}(\lambda_{j}+\gamma)^{s} t^{\alpha}\right)\left(u_{0}, \varphi_{j}\right) \varphi_{j},
	\end{align}
	and the solution operator $\bar{E}$ for $\chi \in L^{2}(\Omega)$ is defined by
	\begin{align}\notag
		\bar{E}(t) \chi=\sum_{j=1}^{\infty} t^{\alpha-1} E_{\alpha, \alpha}\left(-\kappa_{s}(\lambda_{j}+\gamma)^{s} t^{\alpha}\right)\left(\chi, \varphi_{j}\right) \varphi_{j}.
	\end{align}
	
	Before studying the error analysis of full discrete scheme, we also need to give two important properties on the function $E_{\alpha,\beta}$, which can be seen in \cite[Theorem 1.6, formula 1.99]{Poi1999FD}.
	\begin{lemma}\label{lem:3.1}
		Let $0<\alpha<1$, $\beta \in \mathbb{R}$, and $\frac{\alpha \pi}{2}<\mu<\pi\alpha$. Then there exists constant $C=C\left(\alpha, \beta, \mu\right)>0$ such that
		\begin{align}\notag
			E_{\alpha, \beta}\left(z\right) \leq \frac{C}{1+\left|z\right|}, \quad \mu \leq\left|\arg \left(z\right)\right| \leq \pi,~z\in \mathbb{C}.
		\end{align}
	\end{lemma}
	\begin{lemma}\label{lem:3.2}
		Suppose that $\eta >0$, $\mu>0$, and $\alpha>0$. Then one has
		\begin{align}\notag
			\int_{0}^{\eta}\tau^{\alpha-1} E_{\alpha, \alpha}\left(-\mu\tau^{\alpha}\right)d \tau=-\frac{1}{\mu}\int_{0}^{\eta}\frac{d}{d\tau}E_{\alpha, 1}\left(-\mu\tau^{\alpha}\right)d \tau=\frac{1}{\mu}\left(1-E_{\alpha, 1}\left(-\mu\eta^{\alpha}\right)\right).
		\end{align}
	\end{lemma}
	
	It is easy to check that $(-\Delta+\gamma\mathbb{I})$ is a positive self-adjoint operator. Then for all $u \in\dot{H}^{2s}(\Omega)$ and $v$ belonging to a Hilbert space, using Balakrishnan's integral representation, we denote
	\begin{align}\label{equ:3.12}
		\mathcal{A}(u,v):=\left((-\Delta+\gamma\mathbb{I})^{s} u, v\right)=\frac{\sin (\pi s)}{\pi} \int_0^{\infty}\left(z^{s-1}(-\Delta+\gamma\mathbb{I})(z \mathbb{I}+(-\Delta+\gamma\mathbb{I}))^{-1} u, v\right) dz.
	\end{align}
	From \cite{Mag2020FE}, introducing a elliptic projection $P_{h}: H_{0}^{r+1}(\Omega) \rightarrow V_{h}^{d}$ such that for $u\in H_{0}^{r+1}(\Omega)$, we have
	\begin{align}\notag
	\mathcal{A}(P_{h}u-u,v_{h})=0, \quad \forall v_{h} \in V_{h}^{d}.
	\end{align}
	Hence, the following projection estimate holds
	\begin{align}\label{equ:3.13}
	\left\|v-P_{h} v\right\| \leq Ch^{r+1}\|v\|_{r+1}, \quad \forall v \in H_{0}^{r+1}(\Omega),
	\end{align}
	where $h$ denotes maximum of $h_{j}$ for $j=1,2,...,d$. Besides, we also define a linear bijection $\pi_{h}:V_{h}^{d}\to \mathbb{R}^{m}$ such that $\pi_{h}(u_{h})=\bm{u}$, where $\bm{u}$ means solution vector, and $m=\prod\limits_{j=1}^{d}\left(N_{j}-1\right)$ corresponds to the number of degrees of freedom.
	
	Though the proposed fast finite element scheme is based on the linear interpolation, the following discussions can also be extended to high order polynomial interpolation. Next, we primarily study the finite element error analysis of full discrete scheme, but it is also valid when $r\in \mathbb{N}^{+}$ and $r\geq 1$.  
	\begin{lemma}\label{lem:3.3} Let $(-\Delta_{h}+\gamma \mathbb{I})^{s}=\pi_{h}^{-1} (M^{-1}S+\gamma I)^{s} \pi_{h}$ denote discrete spectral fractional Laplacian operator. Then the following estimation holds:
		\begin{align}\notag
		\left| \mathcal{A}(\varphi_{j},v_{h})-\mathcal{A}_{h}(P_{h}\varphi_{j},v_{h})\right| \leq C h^{r+1}\left(\lambda_{j}+\gamma\right)^{s}\left\|\varphi_{j}\right\|_{r+1}\| v_{h} \|,\quad \forall v_{h}\in V_{h}^{d},
		\end{align}
	where $C$ is a positive constant independent of $h$, and $\mathcal{A}_{h}(\cdot,\cdot)$ is defined by
	\begin{align}\notag
		\mathcal{A}_{h}(u_{h},v_{h}):=\left((-\Delta_{h}+\gamma \mathbb{I})^{s}u_{h},v_{h}\right),\quad u_{h},v_{h}\in V_{h}^{d}.
	\end{align}
	\end{lemma}
	
	It suffices to prove Lemma \ref{lem:3.3} by (\ref{equ:3.12}), (\ref{equ:3.13}) and a similar technique developed in \cite[Proposition 1]{Mag2020FE}.
	\begin{lemma}\label{lem:3.4} Let $\left\{\rho\left(t_{k}\right)\right\}_{k=0}^{n}$ denote a series for $n=1,2,...,N$, and define 
		\begin{align}\notag
		\sum_{k=0}^{n}b_{k}^{n}\rho(t_{k}):=\frac{\rho\left(t_{n}\right)}{\tau_{n}}-\frac{\alpha \rho\left(t_{n-1}\right)}{\tau_{n}}-\frac{(1-\alpha)}{\tau_{n}}\left(\frac{\Delta t_{n}}{t_{n}}\right)^{\alpha} \rho\left(t_{0}\right)+\sum_{j=1}^Q w_{j} e^{\xi_{j} \Delta t_{n}}Y_{j}(t_{n-1}).
		\end{align}
		Then we have $b_{n}^{n}=\frac{1}{\tau_{n}}$. In addition, if $n=1$, it follows that
		\begin{align}\notag
		b_{0}^{1}=-\frac{\alpha}{\tau_{1}}-\frac{(1-\alpha)}{\tau_{1}}\left(\frac{\Delta t_{1}}{t_{1}}\right)^{\alpha}=-\frac{1}{\tau_{1}}.
		\end{align}
		If $n=2$, there holds
		\begin{align}\notag
			b_{0}^{2}=\sum\limits_{j=1}^{Q} w_{j} e^{\xi_{j} \Delta t_{2}} k_{2}-\frac{(1-\alpha)}{\tau_{2}}\left(\frac{\Delta t_{2}}{t_{2}}\right)^{\alpha}, \quad
			b_{1}^{2}=\sum\limits_{j=1}^{Q} w_{j} e^{\xi_{j} \Delta t_{2}} k_{3}-\frac{\alpha}{\tau_{2}}.
		\end{align}
		If $n \geq 3$, one has
		\begin{align}\notag
				&b_{0}^{n}=\sum\limits_{j=1}^{Q} w_{j} e^{\xi_{j} \Delta t_{n}} k_{1}^{n-2} k_{2}-\frac{(1-\alpha)}{\tau_{n}}\left(\frac{\Delta t_{n}}{t_{n}}\right)^{\alpha},\quad
				b_{n-1}^{n}=\sum\limits_{j=1}^{Q} w_{j} e^{\xi_{j} \Delta t_{n}} k_{3}-\frac{\alpha}{\tau_{n}}, \\
				&b_{i}^{n}=\sum\limits_{j=1}^{Q} w_{j} e^{\xi_{j} \Delta t_{n}}\left(k_{1} k_{3}+k_{2}\right) k_{1}^{n-2-i}, \quad i=1,2, \ldots, n-2.\notag
		\end{align}
	\end{lemma}
	
	To be more convenient for our analysis, we shall denote
	\begin{align}\label{equ:3.14}
		\bar{c}_{n,1}:=-b_{0}^{n};\quad \bar{c}_{n,n}:=b_{n}^{n};\quad \bar{c}_{n,k}:=\bar{c}_{n,k-1}-b_{k-1}^{n},\quad k=2,3,...,n.
	\end{align}
	Then it can be noted that $\left\{\bar{c}_{n,k}\right\}_{k=1}^{n}$ enjoy explicit representation using $\left\{b_{k}^{n}\right\}_{k=1}^{n}$ in Lemma \ref{lem:3.4}. In addition, it also follows that 
	\begin{align}\notag
		\bar{c}_{n,k}\geq \bar{c}_{n,k-1}\geq \cdots \geq \bar{c}_{n,1}\geq \frac{t_{n}^{-\alpha}}{\Gamma(1-\alpha)}>0,\quad k=1,2,...,n.
	\end{align}
	To prove the above inequality, it just needs to show that $b_{k-1}^{n}\leq 0$ for $k=1,2,...,n$. According to monotonicities of $k_{1}$, $k_{2}$ and $k_{3}$, we have $k_{1}\geq 0$, $k_{2}\geq 0$, $k_{3}\geq 0$, $k_{1}^{n-2}k_{2}\geq 0$ and $\left(k_{1}k_{3}+k_{2}\right)k_{1}^{n-2-i}\geq 0$ in Lemma \ref{lem:3.4}. Then combining negativity of $\left\{w_{j}\right\}_{j=1}^{Q}$ in Lemma \ref{lem:3.4} and (\ref{equ:3.3}), we can obtain the desired result.
	 
	Now, we present the full discrete scheme based on fast time-stepping L1 scheme in time and conforming FEM in space. Find $u_{h}^{n}\in V_{h}^{d}$ such that
	\begin{align}\label{equ:3.15}
		\left(\bar{\partial}_{t}^{\alpha}u_{h}^{n}, v_{h}\right)+\kappa_{s}\mathcal{A}_{h}\left(u_{h}^{n}, v_{h}\right)=\left(f, v_{h}\right),\quad \forall v_{h}\in V_{h}^{d},
    \end{align}
	where $\mathcal{A}_{h}\left(\cdot,\cdot\right)$ is given in Lemma \ref{lem:3.3}, and discrete differential operator $\bar{\partial}_{t}^{\alpha}$ in (\ref{equ:3.7}) is rewritten as
	\begin{align}\label{equ:3.16}
		\bar{\partial}_{t}^{\alpha}u_{h}^{n}=\sum_{k=1}^{n}\bar{c}_{n,k}\left(u_{h}^{k}-u_{h}^{k-1}\right),\quad n=1,2,...,N.
	\end{align}
	
	Notice that full discrete scheme (\ref{equ:3.15}) does not require $P_{h}u(t_{0})=u_{h}^{0}$. Below, we first present a lemma to derive stability analysis.
	\begin{lemma}\label{lem:3.5} Let $\left\{v^{j}\right\}_{j=0}^{m}$ denote arbitrary series. Then for $m=1,2,...,n$, one has
		\begin{align}\notag
			&\mathcal{A}_{h}(v^{m},\sum_{k=1}^{m}\bar{c}_{n,k}\left(v^{k}-v^{k-1}\right))\geq\frac{1}{2} \sum_{k=1}^{m}\bar{c}_{n,k}\left(\mathcal{A}_{h}(v^{k},v^{k})-\mathcal{A}_{h}(v^{k-1},v^{k-1})\right),\\
			&(v^{m},\sum_{k=1}^{m}\bar{c}_{n,k}\left(v^{k}-v^{k-1}\right))\geq\Vert v^{m}\Vert\sum_{k=1}^{m}\bar{c}_{n,k}\left(\Vert v^{k}\Vert-\Vert v^{k-1}\Vert\right),\notag
		\end{align}
	where $\mathcal{A}_{h}(\cdot,\cdot)$ is defined in Lemma \ref{lem:3.3}, and $n=1,2,...,N$.
	\end{lemma}
		\begin{proof}
		We temporarily let
		\begin{align}\notag
			H_{m}=\mathcal{A}_{h}(v^{m},\sum_{k=1}^{m}\bar{c}_{n,k}\left(v^{k}-v^{k-1}\right))-\frac{1}{2} \sum_{k=1}^{m}\bar{c}_{n,k}\left(\mathcal{A}_{h}(v^{k},v^{k})-\mathcal{A}_{h}(v^{k-1},v^{k-1})\right).
		\end{align}
		Then it just needs to show that $H_{m}\geq 0$. Since case $m=1$ is clear, we discuss $m\geq 2$, namely,
		\begin{align}\label{equ:3.17}
			H_{m}&=\sum_{k=1}^{m}\bar{c}_{n,k}\mathcal{A}_{h}\left(v^{k}-v^{k-1},\frac{1}{2}\left(v^{k}-v^{k-1}\right)+\sum_{s=k+1}^{m}(v^{s}-v^{s-1})\right) \notag\\
			&=\frac{1}{2}\sum_{k=1}^{m}\bar{c}_{n,k}\mathcal{A}_{h}\left(v^{k}-v^{k-1},v^{k}-v^{k-1}\right)+
			\sum_{k=1}^{m-1}\bar{c}_{n,k}\mathcal{A}_{h}\left(v^{k}-v^{k-1},\sum_{s=k+1}^{m}(v^{s}-v^{s-1})\right) \notag\\
			&=\frac{1}{2}\sum_{k=1}^{m}\bar{c}_{n,k}\mathcal{A}_{h}\left(v^{k}-v^{k-1},v^{k}-v^{k-1}\right)+
			\sum_{s=1}^{m}\mathcal{A}_{h}\left(v^{s}-v^{s-1},\sum_{k=1}^{s-1}\bar{c}_{n,k}(v^{k}-v^{k-1})\right),
		\end{align}
		where the last equality used the following fact 
		\begin{align}\notag
			\sum_{k=1}^{m-1}\mathcal{A}_{h}\left(z^{k},\sum_{s=k+1}^{m}w^{s}\right)=\sum_{s=1}^{m}\mathcal{A}_{h}\left(w^{s},\sum_{k=1}^{s-1}z^{k}\right).
		\end{align}
		Now we let $z^{0}=0$, and $z^{k}-z^{k-1}=\left(v^{k}-v^{k-1}\right)\bar{c}_{n,k}$ for $k=1,2,...,n$. From (\ref{equ:3.17}) we derive 
		\begin{align}\notag
			H_{m}&=\frac{1}{2}\sum_{k=1}^{m}\frac{1}{\bar{c}_{n,k}}\mathcal{A}_{h}\left(z^{k}-z^{k-1},z^{k}-z^{k-1}\right)+
			\sum_{s=1}^{m}\frac{1}{\bar{c}_{n,s}}\mathcal{A}_{h}\left(z^{s}-z^{s-1},z^{s-1}\right) \notag\\
			&=\frac{1}{2}\sum_{k=1}^{m-1}\left(\frac{1}{\bar{c}_{n,k}}-\frac{1}{\bar{c}_{n,k+1}}\right)\mathcal{A}_{h}\left(z^{k},z^{k}\right)+\frac{1}{2}\frac{1}{\bar{c}_{n,m}}\mathcal{A}_{h}\left(z^{m},z^{m}\right)\geq 0.\notag
		\end{align}
		However, the second conclusion can be obtained by a direct expansion and the use of Cauchy-Schwarz inequality. Hence, we complete this proof of the lemma.
	\end{proof}
	
	\begin{theorem}\label{the:3.1} Assume that $u_{h}^{n}\in V_{h}^{d}$ corresponds to the solution of full discrete scheme (\ref{equ:3.15}) using (\ref{equ:3.16}). Then we have
		\begin{align}\notag
			\Vert u_{h}^{n}\Vert \leq \frac{1}{\bar{c}_{n,n}}\left(\Vert f^{n}\Vert+\bar{c}_{n,1}\Vert u_{h}^{0}\Vert+\sum_{k=1}^{n-1}(\bar{c}_{n,k+1}-\bar{c}_{n,k})\Vert u_{h}^{k}\Vert\right),\quad n=1,2,...,N,
		\end{align}
		where $1/\bar{c}_{n,n}=\tau_{n}$.
	\end{theorem}
	\begin{proof}
		From (\ref{equ:3.15}), we take $v_{h}=u_{h}^{n}$ and invoke Lemma \ref{lem:3.5}, then it follows that
		\begin{align}\notag
			\Vert u_{h}^{n}\Vert\sum_{k=1}^{n}\bar{c}_{n,k}\left(\Vert u_{h}^{k}\Vert-\Vert u_{h}^{k-1}\Vert\right)+\kappa_{s}\mathcal{A}_{h}\left(u_{h}^{n},u_{h}^{n}\right)\leq \left(\bar{\partial}_{t}^{\alpha}u_{h}^{n},u_{h}^{n}\right)
			+\kappa_{s}\mathcal{A}_{h}\left(u_{h}^{n},u_{h}^{n}\right)=\left(f^{n},u_{h}^{n}\right).
		\end{align}
	The uses of Cauchy-Schwarz inequality and positivity of $\mathcal{A}_{h}(\cdot,\cdot)$ give
	\begin{align}\notag
		\sum_{k=1}^{n}\bar{c}_{n,k}\left(\Vert u_{h}^{k}\Vert-\Vert u_{h}^{k-1}\Vert\right)\leq \Vert f^{n}\Vert,\quad n=1,2,...,N.
	\end{align}
	After a simple recombination of the above inequality yields the $L^{2}$ stability estimate.
	\end{proof}
	
	Below is an $L^{2}$ stability with respect to a weighted sum. Then we define
	\begin{align}\notag
		\theta_{n,n}=1, \quad and\quad \theta_{n,j}=\Gamma(2-\alpha)\sum\limits_{k=j}^{n-1}\Delta t_{k}^{\alpha}\theta_{k,j}\left(\bar{c}_{n,k+1}-\bar{c}_{n,k}\right),\quad j=1,2,...,n-1.
	\end{align}
	\begin{theorem}\label{the:3.2} Suppose that $u_{h}^{n}\in V_{h}^{d}$ denotes the solution of full discrete scheme (\ref{equ:3.15}). Then we have
		\begin{align}\notag
			\Vert u_{h}^{n}\Vert\leq\Vert u_{h}^{0}\Vert+\frac{1}{\bar{c}_{n,n}}\sum_{j=1}^{n}\theta_{n,j}\Vert f^{j}\Vert,
		\end{align}
	  where $1/\bar{c}_{n,n}=\tau_{n}$, and $n=1,2,...,N$.
	\end{theorem}
	\begin{proof}
		To complete this proof, we now apply mathematical induction on $k$, $k=1,2,...,n$. When $k=1$ it is clear. Assume that conclusion holds for $k=1,2,...,n-1$. Then invoking Theorem \ref{the:3.1}, we shall obtain 
		\begin{align}
			\Vert u_{h}^{n}\Vert&\leq \frac{1}{\bar{c}_{n,n}}\left(\Vert f^{n}\Vert+\bar{c}_{n,1}\Vert u_{h}^{0}\Vert+\sum_{k=1}^{n-1}(\bar{c}_{n,k+1}-\bar{c}_{n,k})\Vert u_{h}^{k}\Vert\right) \notag\\
			&\leq \frac{1}{\bar{c}_{n,n}}\left(\Vert f^{n}\Vert+\bar{c}_{n,n}\Vert u_{h}^{0}\Vert+\Gamma(2-\alpha)\sum_{j=1}^{n-1}\sum_{k=j}^{n-1}(\bar{c}_{n,k+1}-\bar{c}_{n,k})\theta_{k,j}\Delta t_{k}^{\alpha}\Vert f^{j}\Vert\right) \notag\\
			&= \frac{1}{\bar{c}_{n,n}}\left(\Vert f^{n}\Vert+\bar{c}_{n,n}\Vert u_{h}^{0}\Vert+\sum_{j=1}^{n-1}\theta_{n,j}\Vert f^{j}\Vert\right),\notag
		\end{align}
		where the second line is obtained using equality
		\begin{align}\notag
			\sum\limits_{k=1}^{n-1}w_{n,k}\sum\limits_{j=1}^{k}z_{k,j}=\sum\limits_{j=1}^{n-1}\sum\limits_{k=j}^{n-1}w_{n,k}z_{k,j}.
		\end{align}
		Therefore, we get the desired result.
	\end{proof}
	
	Theorem \ref{the:3.2} shows that the stability bound shall not blow up as $\alpha \to 1^{-}$. To obtain the error estimate for full discrete scheme based on fast time-stepping scheme on graded mesh, we also need to provide an useful lemma. Though the definition of $\theta_{n,j}$ is different from direct L1 case in \cite[Lemma 5.3, Corollary 5.4, Corollary 5.5]{Chh2021BU}, we can derive the following results.
	\begin{lemma}\label{lem:3.6} Assume $1\leq\omega \leq 2(2-\alpha)/\alpha$, and let $n=1,2, \ldots, N$. (a) Then for $\eta \in(0,1)$, we can have
	\begin{align}\notag
		\Delta t_{n}^\alpha \sum_{j=1}^n j^{\omega(\eta-\alpha)} \theta_{n, j} \leq \frac{\Gamma(1-\alpha)}{\Gamma(2-\alpha)} \left(\frac{t_n}{T}\right)^\eta N^{\omega(\eta-\alpha)}T^{\alpha}.
	\end{align}
	(b) When taking $\eta=\alpha$ in (a), one has
	\begin{align}\notag
		\Delta t_{n}^\alpha \sum_{j=1}^n\theta_{n, j} \leq \frac{\Gamma(1-\alpha)}{\Gamma(2-\alpha)}t_{n}^{\alpha}.
	\end{align}
	(c) Denote $m^{*}=\left\{2-\alpha,\omega\alpha\right\}$, and choose $\eta=\delta+\alpha-\frac{m^*}{\omega}$ for $\delta\in\left(0,1-\frac{\alpha}{2}\right)$. There holds
	\begin{align}\notag
	\Delta t_{n}^\alpha \sum_{j=1}^n j^{-m^*} \theta_{n, j} \leq \frac{\Gamma(1-\alpha)}{\Gamma(2-\alpha)}N^{\omega\delta} T^\alpha\left(\frac{t_n}{T}\right)^{\delta+\alpha-\frac{m^*}{\omega}} N^{-m^*}\notag.
	\end{align}
	\begin{proof}
		We now start to prove (a). Since $\bar{c}_{j,k}\geq\frac{t_{j}^{-\alpha}}{\Gamma(1-\alpha)}$ for $k=1,2,...,j$, we can derive
		\begin{align}
			\bar{\partial}_{t}^{\alpha}t_{j}^{\eta}=\sum_{k=1}^{j}\bar{c}_{j,k}\left(t_{k}^{\eta}-t_{k-1}^{\eta}\right)
			=\eta\sum_{k=1}^{j}\bar{c}_{j,k}\int_{t_{k-1}}^{t_{k}}\zeta^{\eta-1}d\zeta\geq\frac{t_{j}^{\eta-\alpha}}{\Gamma(1-\alpha)},\quad j=1,2,...,n.\notag
		\end{align}
		Multiplying the above inequality by $\theta_{n,j}$, and summing from $j=1$ to $n$ yield
		\begin{align}
			\frac{1}{\Gamma(1-\alpha)}\sum_{j=1}^{n}t_{j}^{\eta-\alpha}\theta_{n,j}&\leq
			\sum_{j=1}^{n}\theta_{n,j}\sum_{k=1}^{j}\bar{c}_{j,k}\left(t_{k}^{\eta}-t_{k-1}^{\eta}\right)
			=\sum_{k=1}^{n}\left(t_{k}^{\eta}-t_{k-1}^{\eta}\right)\sum_{j=k}^{n}\bar{c}_{j,k}\theta_{n,j}\\
			&=\frac{1}{\Delta t_{n}^{\alpha}\Gamma(2-\alpha)}\sum_{k=1}^{n}\left(t_{k}^{\eta}-t_{k-1}^{\eta}\right)
			=\frac{t_{n}^{\eta}}{\Delta t_{n}^{\alpha}\Gamma(2-\alpha)},\notag
		\end{align}
		where we used mathematical induction to prove
		\begin{align}
			\frac{1}{\Delta t_{n}^{\alpha}\Gamma(2-\alpha)}=\sum_{j=k}^{n}\bar{c}_{j,k}\theta_{n,j},\quad 1\leq k\leq n.\notag
		\end{align}
	    Then we obtain (a) in terms of $t_{j}=\left(\frac{j}{N}\right)^{\omega}T$. The conclusion (b) is clear and the proof of (c) similarly follows from \cite[Corollary 5.5]{Chh2021BU}.
	\end{proof}
	\end{lemma}
	
	To explain why we consider graded time mesh and use the following time regularity assumption in the present temporal convergence analysis, a reason is illustrated in Appendix. Now, we shall give the following time regularity assumption.
	\begin{assumption}\label{ass:3.1}\cite{Chh2021BU,Stm2017EA} Let initial condition $u_{0}$ and right hand side function $f$ in problem (\ref{equ:1.1}) satisfy appropriate conditions such that the solution $u$ in (\ref{equ:3.11}) of (\ref{equ:1.1}) has
		\begin{align}\notag
			\Vert \partial_{t}^{l}u\Vert_{q}\leq C\left(1+t^{\alpha-l}\right),\quad t>0,
		\end{align}
		where $C$ is a positive constant, $l=0,1,2$, and $q$ denotes non-negative integer number.
	\end{assumption}
	
	Analogous to \cite[Lemma 2.6]{Jis2017FE}, using \cite[Lemma 5.2]{Stm2017EA}, Assumption \ref{ass:3.1} and triangle inequality, we can bound
	\begin{align}\label{equ:3.18}
		\Vert r^{n}\Vert_{q}\leq C\left(n^{-\min\left\{2-\alpha,\omega\alpha\right\}}+\epsilon\right),\quad n=1,2,...,N;\quad q=0,1,
	\end{align}
	where $r^{n}$ denotes the truncation error of fast time-stepping L1 scheme in (\ref{equ:3.7}), $\epsilon$ is given in (\ref{equ:3.5}) to satisfy (\ref{equ:3.4}), and $C$ is a positive constant.
	 
	We next extend the technique developed in \cite{Mag2020FE} to time-space fractional diffusion problem with imhomogenous initial value and right hand side term. For simplicity, we temporarily denote $u_{j}(t):=(u(\cdot,t),\varphi_{j})$, $f_{j}(t):=(f(\cdot,t),\varphi_{j})$ and $u_{0j}:=(u_{0}(\cdot),\varphi_{j})$. 
	\begin{theorem}\label{the:3.3} Suppose that $\nu_{0}=\sum\limits_{j=1}^{\infty}\left|u_{0j}\right|\left(\lambda_j+\gamma\right)^{\frac{r+1}{2}+s}<\infty$ and $\nu_{1}=\sup\limits_{0 \leq t \leq T} \sum\limits_{j=1}^{\infty}\left(\lambda_{j}+\gamma\right)^{\frac{r+1}{2}}\left|f_{j}(t)\right|<\infty$. If $u(t)$ in  (\ref{equ:3.11}) is the solution of (\ref{equ:1.1}), then we have
		\begin{align}\notag
			&\left|\left((-\Delta+\gamma\mathbb{I})^{s} u(t)-\left(-\Delta_{h}+\gamma\mathbb{I}\right)^{s} P_{h} u(t), v_{h}\right)\right|\leq C\left(\nu_{0}+\nu_{1}\right) h^{r+1}\left\|v_{h}\right\|,\quad \forall v_{h}\in V_{h}^{d},
		\end{align}
	for $t>0$, where $C$ is a positive constant.
	\end{theorem}
	\begin{proof} By Lemma \ref{lem:3.1} and Lemma \ref{lem:3.2}, we have
		\begin{align}\label{equ:3.19}
			\left|u_{j}(t)\right| &=\left|u_{0j} E_{a,1}\left(-\kappa_{s}\left(\lambda_{j}+\gamma\right)^{s} t^{\alpha}\right)+\int_{0}^{t}(t-\tau)^{\alpha-1} E_{\alpha, \alpha}\left(-\kappa_{s}\left(\lambda_{j}+\gamma\right)^{s}(t-\tau)^{\alpha}\right) f_{j}(\tau)d \tau\right| \notag\\
			& \leq\left|u_{0j}\right| \frac{C_1}{1+\kappa_{s}\left(\lambda_{j}+\gamma\right)^{s} t^{\alpha}}+\sup_{0 \leq \tau \leq T}\left|f_{j}(\tau)\right|\left|\int_{0}^{t}(t-\tau)^{\alpha-1} E_{\alpha, \alpha}\left(-\kappa_{s}\left(\lambda_{j}+\gamma\right)^{s}(t-\tau)^{\alpha}\right) d \tau\right| \notag\\
			&  \leq C_{1}\left|u_{0j}\right|+\frac{C_{2}}{\left(\lambda_{j}+\gamma\right)^{s}} \sup _{0\leq \tau \leq T}\left|f_{j}(\tau)\right|,\quad \forall t>0.
		\end{align}
	In addition, the $H^{r+1}$-norm of $\varphi_{j}$ satisfies the following growth condition:
		\begin{align}\label{equ:3.20}
			\left\|\varphi_{j}\right\|_{r+1}\leq C\lambda_{j}^{\frac{r+1}{2}}\leq C \left(\lambda_{j}+\gamma\right)^{\frac{r+1}{2}},\quad \forall j\in \mathbb{N}^{+}.
		\end{align}
	Hence, using (\ref{equ:3.19}), (\ref{equ:3.20}) and Lemma \ref{lem:3.3}, we can derive the estimate
		\begin{align}
			&\left|\left(\left((-\Delta+\gamma\mathbb{I})^{s}-\left(-\Delta_{h}+\gamma\mathbb{I}\right)^{s} P_{h} \right)u(t),v_{h}\right)\right|
			\leq \sum_{j=1}^{\infty}\left|\left(u_{j}(t)\left((-\Delta+\gamma\mathbb{I})^{s}-\left(-\Delta_{h}+\gamma\mathbb{I}\right)^{s} P_{h}\right) \varphi_{j}, v_{h}\right)\right| \notag\\
			&\leq C h^{r+1}\left\|v_{h}\right\| \sum_{j=1}^{\infty}\left\|\varphi_{j}\right\|_{k}\left|u_{j}(t)\right|\left(\lambda_{j}+\gamma\right)^{s} \leq C h^{r+1}\left\|v_{h}\right\| \sum_{j=1}^{\infty}\left|u_{j}(t)\right|\left(\lambda_{j}+\gamma\right)^{\frac{r+1}{2}+s} \notag\\
			&\leq C h^{r+1}\left\|v_{h}\right\|\left(\nu_{0}+\nu_{1}\right),\quad \forall v_{h}\in V_{h}^{d}.\notag
		\end{align}
	This completes the proof.
	\end{proof}
	\begin{theorem}\label{the:3.4} Let $u(t)$ be the solution of (\ref{equ:1.1}) satisfying regularity in Assumption \ref{ass:3.1}, graded mesh factor $\omega \leq 2(2-\alpha)/\alpha$, and $u(t_{n})\in H^{r+1}\left(\Omega\right)$ for $n=0,1,...,N$. Under the assumptions of Theorem \ref{the:3.3} and Lemma \ref{lem:3.6}, the full discrete scheme (\ref{equ:3.15}) has the following error estimate
		\begin{align}\notag
			\max_{1 \leq n \leq N}\left\|u(t_{n})-u_{h}^{n}\right\|\leq C\left(N^{-\left\{2-\alpha,\omega\alpha\right\}}+h^{r+1}+\left\|u(t_{0})-u_{h}^{0}\right\|+\epsilon\right),
		\end{align}
		where $C$ corresponds to a positive constant, and $\epsilon$ is given in (\ref{equ:3.5}) to satisfy (\ref{equ:3.4}) and it is small sufficiently enough to be ignored.
	\end{theorem}
	\begin{proof} We first derive semi-discrete scheme for simplicity of analysis. Find $u(t_{n})\in \dot{H}^{2s}(\Omega)$ such that
		\begin{align}\label{equ:3.21}
			\left(\partial_{t}^{\alpha}u(t_{n}), v_{h}\right)+\kappa_{s}\left((-\Delta+\gamma \mathbb{I})^{s} u(t_{n}), v_{h}\right)=\left(f^{n}, v_{h}\right),\quad \forall v_{h}\in V_{h}^{d},
		\end{align}
	    where $f^{n}=f(\cdot,t_{n})$. According to (\ref{equ:3.21}), we immediately have
		\begin{align}\label{equ:3.22}
			\left(\bar{\partial}_{t}^{\alpha}P_{h} u(t_{n}), v_{h}\right)&+\kappa_{s}\left(\left(-\Delta_{h}+\gamma \mathbb{I}\right)^{s} P_{h} u(t_{n}), v_{h}\right) \notag\\
			&=\left(\left(\bar{\partial}_{t}^{\alpha}-\partial_{t}^{\alpha}\right)u(t_{n}), v_{h}\right)+\left(f^{n}, v_{h}\right)
			+\left(\bar{\partial}_{t}^{\alpha}\left(P_{h}-\mathbb{I}\right)u(t_{n}), v_{h}\right) \notag\\
			&+\kappa_{s}\left(\left(\left(-\Delta_{h}+\gamma \mathbb{I}\right)^{s} P_{h}-(-\Delta+\gamma \mathbb{I})^{s}\right)u(t_{n}), v_{h}\right),\quad \forall v_{h}\in V_{h}^{d}.	
		\end{align}
	Combining (\ref{equ:3.16}), (\ref{equ:3.13}) and Assumption \ref{ass:3.1} with $q=r+1$ and $l=1$ easily gives
	\begin{align}\label{equ:3.23}
		\left\|\bar{\partial}_{t}^{\alpha}\left(P_{h}-\mathbb{I}\right)u(t_{n})\right\|
		&\leq Ch^{r+1}\sum_{k=1}^{n}\bar{c}_{n,k}\int_{t_{k-1}}^{t_{k}}\Vert u_{\tau}\Vert_{r+1} d\tau
		\leq Ch^{r+1}\sum_{k=1}^{n}\bar{c}_{n,k}\int_{t_{k-1}}^{t_{k}}\tau^{\alpha-1} d\tau,\notag\\
		&\leq Ch^{r+1}\bar{c}_{n,n}\int_{0}^{t_{n}}\tau^{\alpha-1}d\tau=Ch^{r+1}.
	\end{align}
	Now, denoting $\rho^{n}=P_{h}u(t_{n})-u_{h}^{n}$, and subtracting (\ref{equ:3.15}) from (\ref{equ:3.22}), we have
	\begin{align}\label{equ:3.24}
		\left(\bar{\partial}_{t}^{\alpha}\rho^{n}, v_{h}\right)&+\kappa_{s}\left(\left(-\Delta_{h}+\gamma \mathbb{I}\right)^{s} \rho^{n}, v_{h}\right) 
		=\left(\left(\bar{\partial}_{t}^{\alpha}-\partial_{t}^{\alpha}\right)u(t_{n}), v_{h}\right)+
		\left(\bar{\partial}_{t}^{\alpha}\left(P_{h}-\mathbb{I}\right)u(t_{n}), v_{h}\right) \notag\\
		&+\kappa_{s}\left(\left(\left(-\Delta_{h}+\gamma \mathbb{I}\right)^{s} P_{h}-(-\Delta+\gamma \mathbb{I})^{s}\right)u(t_{n}), v_{h}\right),\quad \forall v_{h}\in V_{h}^{d}.	
	\end{align}
	
	However, due to the use of Theorem \ref{the:3.3}, we can not directly invoke stability estimate established in Theorem \ref{the:3.2} to bound $\rho^{n}$. Hence, we have to return back to the proof of Theorem \ref{the:3.1} by considering (\ref{equ:3.24}).
	
	Similar to Theorem \ref{the:3.1}, taking $v_{h}=\rho^{n}$ in (\ref{equ:3.24}) and invoking Theorem \ref{the:3.3}, we have 
	\begin{align}\label{equ:3.25}
		\Vert \rho^{n}\Vert\sum_{k=1}^{n}\bar{c}_{n,k}\left(\Vert \rho^{k}\Vert-\Vert \rho^{k-1}\Vert\right)&\leq \Vert \left(\bar{\partial}_{t}^{\alpha}-\partial_{t}^{\alpha}\right)u(t_{n})\Vert\left\|\rho^{n}\right\|+\Vert \bar{\partial}_{t}^{\alpha}\left(P_{h}-\mathbb{I}\right)u(t_{n})\Vert\left\|\rho^{n}\right\|\notag\\
		&+C\kappa_{s}\left(\nu_{0}+\nu_{1}\right)h^{r+1}\left\|\rho^{n}\right\|=:\left\|z^{n}\right\|\left\|\rho^{n}\right\|,
	\end{align}
	where $C$, $\nu_{0}$ and $\nu_{1}$ are three constants given by Theorem \ref{the:3.3}. Hence, eliminating $\left\|\rho^{n}\right\|$ in (\ref{equ:3.25}) yields 
	\begin{align}\notag
		\Vert \rho^{n}\Vert \leq \frac{1}{\bar{c}_{n,n}}\left(\Vert z^{n}\Vert+\bar{c}_{n,1}\Vert \rho^{0}\Vert+\sum_{k=1}^{n-1}(\bar{c}_{n,k+1}-\bar{c}_{n,k})\Vert \rho^{k}\Vert\right),\quad n=1,2,...,N.
	\end{align}
	A same technique in Theorem \ref{the:3.2} derives
	\begin{align}\label{equ:3.26}
		\Vert \rho^{n}\Vert\leq\Vert \rho^{0}\Vert+\Delta t_{n}^{\alpha}\Gamma(2-\alpha)\sum_{j=1}^{n}\theta_{n,j}\Vert z^{j}\Vert,\quad n=1,2,...,N.
	\end{align}
	Next, we are ready to bound $\Vert \rho^{n}\Vert$. Substituting $\Vert z^{j}\Vert$ in (\ref{equ:3.25}) into (\ref{equ:3.26}) leads to
	\begin{align}\label{equ:3.27}
		\Vert \rho^{n}\Vert&\leq\Vert \rho^{0}\Vert+\frac{1}{\bar{c}_{n,n}}\sum_{j=1}^{n}\theta_{n,j}\Vert \left(\bar{\partial}_{t}^{\alpha}-\partial_{t}^{\alpha}\right)u(t_{j})\Vert
		+\frac{1}{\bar{c}_{n,n}}\sum_{j=1}^{n}\theta_{n,j}\Vert \bar{\partial}_{t}^{\alpha}\left(P_{h}-\mathbb{I}\right)u(t_{j})\Vert
		\notag\\
		&+C\kappa_{s}\left(\nu_{0}+\nu_{1}\right)h^{r+1}\frac{1}{\bar{c}_{n,n}}\sum_{j=1}^{n}\theta_{n,j}\notag\\
		&\leq\Vert \rho^{0}\Vert+C\Delta t_{n}^{\alpha}\sum_{j=1}^{n}\theta_{n,j}j^{-\left\{2-\alpha,\omega\alpha\right\}}+
		C\left(h^{r+1}+\epsilon\right)\Delta t_{n}^{\alpha}\sum_{j=1}^{n}\theta_{n,j}\notag\\
		&\leq \Vert\rho^{0}\Vert+C\left(N^{-\left\{2-\alpha,\omega\alpha\right\}}+h^{r+1}+\epsilon\right),\quad n=1,2,...,N,
	\end{align}
	where we also used (\ref{equ:3.18}), (\ref{equ:3.23}) and Lemma \ref{lem:3.6}. By triangle inequality we will also obtain
	\begin{align}\label{equ:3.28}
		\left\|\rho^{0}\right\|\leq \left\|P_{h}u(t_{0})-u(t_{0})\right\|+\left\|u(t_{0})-u_{h}^{0}\right\|\leq Ch^{r+1}\left\|u(t_{0})\right\|_{r+1}+\left\|u(t_{0})-u_{h}^{0}\right\|.
	\end{align}
	Hence, combining (\ref{equ:3.13}), (\ref{equ:3.27}) and (\ref{equ:3.28}), and using
	\begin{align}\notag
		\max_{1 \leq n \leq N}\left\|u(t_{n})-u_{h}^{n}\right\|\leq \max_{1 \leq n \leq N}\left\|u(t_{n})-P_{h}u(t_{n})\right\| 
		+\max_{1 \leq n \leq N}\left\|\rho^{n}\right\|,
	\end{align}
	we immediately prove the desired result.
	\end{proof}
	
	From Theorem \ref{the:3.4}, we can notice that the choice of $\omega=(2-\alpha)/\alpha$ shall yield an optimal temporal convergence $\mathcal{O}(N^{-(2-\alpha)})$, but the error bound is $\alpha$-nonrobust as $\alpha\to 1^{-}$ due to the use of Lemma \ref{lem:3.6}. For the convergence analysis of full discrete scheme (\ref{equ:3.9}) employing fourth-order CDM, a similar discussion can be seen in \cite{Seb2017FD}. Next, we shall provide some numerical examples to illustrate our theoretical achievements.
\section{Numerical experiments}
	In this section, several numerical examples are delivered to test the performance of the present fast algorithm for fractional diffusion equations with spectral fractional Laplacian. Also, all numerical simulations are conducted using Python 3.7 on the personal computer.
	 \begin{table*}[htb]\small
	\setlength{\abovecaptionskip}{0cm}
	\caption{\label{tab:4.1}The $L^{2}$ error and the corresponding convergence rate of Example \ref{exa:4.1} in linear FEM.}
	\centering
	\renewcommand\arraystretch{1.1}
	\begin{tabular}{p{1cm} p{1cm} p{1.7cm} p{1.7cm}p{1.7cm} p{1.7cm} p{1.7cm} p{1.7cm}}
		\hline
		$(s,\gamma)$&                 $\mathcal{M}$              &   8          & 16        & 32        & 64          & 128 & 256 \\
		\hline 
		$(0.4,0)$  & $\Vert e_{h}\Vert_{2}$   &   1.746E-02  & 3.908E-03 & 9.507E-04 & 2.361E-04   & 5.892E-05  & 1.472E-05          \\
		&  $~~r_{2}$               &              &  2.160    &  2.039    & 2.010       & 2.002 & 2.001 \\
		$(1.2,0)$  & $\Vert e_{h}\Vert_{2}$   &   3.218E-04  & 7.312E-05 & 1.786E-05 & 4.438E-06   & 1.108E-06 & 2.769E-07 \\
		& $~~r_{2}$                &              &  2.138    &  2.034    & 2.008       & 2.002 & 2.001 \\
		\hline
		$(0.8,1)$  & $\Vert e_{h}\Vert_{2}$   &   2.365E-03  & 5.332E-04 & 1.300E-04 & 3.229E-05   & 8.059E-06 & 2.014E-06 \\
		& $~~r_{2}$                &              &  2.149    & 2.037     & 2.009       & 2.002  & 2.001\\
		$(1.6,1)$  & $\Vert e_{h}\Vert_{2}$   &   4.282E-05  & 9.802E-06 & 2.398E-06 & 5.962E-07   & 1.489E-07  & 3.720E-08\\
		& $~~r_{2}$                &              &  2.127    &  2.031    & 2.008       & 2.002 & 2.001\\
		\hline
		$(1.0,2)$  & $\Vert e_{h}\Vert_{2}$   &   8.639E-04  & 1.955E-04 & 4.770E-05 & 1.185E-05   & 2.959E-06 & 7.394E-07\\
		& $~~r_{2}$                &              &  2.144    &  2.035    & 2.009       & 2.002 & 2.001\\
		$(2.0,2)$  & $\Vert e_{h}\Vert_{2}$   &   5.628E-06  & 1.297E-06 & 3.179E-07 & 7.909E-08   & 1.975E-08 & 4.936E-09 \\
		& $~~r_{2}$                &              &  2.117    & 2.029     & 2.007       & 2.002 & 2.001 \\
		\hline
	\end{tabular}
    \end{table*}
\subsection{Accuracy test}
	To check the accuracy and efficiency of the present algorithm, we first study the high dimensional fractional Possion problem (\ref{equ:2.1}). Unless otherwise specified, $\Omega=(0,1)^{d}$, $\mathcal{M}=N_{1}=\cdots=N_{d}$, $h=\mathcal{M}^{-1}$ and $\bm{x}=(x_{1},x_{2},...,x_{d})$ are given throughout numerical section.
	\begin{table*}[htb]\small
		\setlength{\abovecaptionskip}{0cm}
		\caption{\label{tab:4.2} $L^{2}$ errors and the corresponding convergence rates of Example \ref{exa:4.1} in fourth-order CDM.}
		\centering
		\renewcommand\arraystretch{1.1}
		\begin{tabular}{p{1cm} p{1cm} p{1.7cm} p{1.7cm}p{1.7cm} p{1.7cm} p{1.7cm} p{1.7cm}}
			\hline
			$(s,\gamma)$&                 $\mathcal{M}$              &   8          & 16        & 32        & 64          & 128 & 256 \\
			\hline 
		    $(0.3,0)$  & $\Vert e_{h}\Vert_{2}$   &   4.113E-05  & 2.525E-06 & 1.571E-07 & 9.806E-09   & 6.127E-10  & 3.829E-11          \\
			&  $~~r_{2}$               &              &  4.026    &  4.007    & 4.002       & 4.000 & 4.000 \\
			$(1.3,0)$  & $\Vert e_{h}\Vert_{2}$   &   1.506E-06  & 9.238E-08 & 5.747E-09 & 3.588E-10   & 2.242E-11 & 1.401E-12 \\
			& $~~r_{2}$                &              &  4.027    &  4.007    & 4.002       & 4.000 & 4.000 \\
			\hline
			$(0.5,1)$  & $\Vert e_{h}\Vert_{2}$   &   2.606E-05  & 1.599E-06 & 9.949E-08 & 6.211E-09   & 3.881E-10 & 2.425E-11 \\
			& $~~r_{2}$                &              &  4.026    & 4.007     & 4.002       & 4.000  & 4.000\\
			$(1.5,1)$  & $\Vert e_{h}\Vert_{2}$   &   6.550E-07  & 4.017E-08 & 2.499E-09 & 1.560E-10   & 9.748E-12  & 6.092E-13           \\
			& $~~r_{2}$                &              &  4.027    & 4.007     & 4.002       & 4.000  & 4.000\\
			\hline
			$(0.7,2)$  & $\Vert e_{h}\Vert_{2}$   &   1.382E-05  & 8.481E-07 & 5.276E-08 & 3.294E-09   & 2.058E-10 & 1.286E-11 \\
			& $~~r_{2}$                &              &  4.026    & 4.007     & 4.002       & 4.000  & 4.000\\
			$(1.7,2)$  & $\Vert e_{h}\Vert_{2}$   &   2.789E-07  & 1.710E-08 & 1.064E-09 & 6.642E-11   & 4.150E-12 & 2.594E-13 \\
			& $~~r_{2}$                &              &  4.028    & 4.007     & 4.001       & 4.000  & 4.000\\
			\hline
		\end{tabular}
	\end{table*}
	
	From the definition of $(-\Delta+\gamma \mathbb{I})^{s}$, it will be convenient to construct the exact solution by replacing $j$ with $(i_{1},i_{2},...,i_{d})$ in (\ref{equ:1.3}). Hence the exact solution of (\ref{equ:2.1}) can be easily expressed as
	\begin{align}\label{equ:4.1}
		u(\bm{x})=\sum_{i_{1},i_{2},...,i_{d}=1}^{\infty} \left(\lambda_{i_{1},i_{2},...,i_{d}}+\gamma\right)^{-s} \widehat{f}_{i_{1},i_{2},...,i_{d}} \varphi_{i_{1},i_{2},...,i_{d}}(\bm{x}), \quad \forall\bm{x} \in \bar{\Omega},
	\end{align}
	where the coefficients
	\begin{align}\notag
		\widehat{f}_{i_{1},i_{2},...,i_{d}}=\int_{\Omega} f(\bm{x}) \varphi_{i_{1},i_{2},...,i_{d}}(\bm{x}) d\bm{x}, \quad i_{1},i_{2},...,i_{d} \in \mathbb{N}^{+} .
	\end{align}
	Besides, we choose the eigenpairs of Laplacian operator $(-\Delta)$ as follows:
	\begin{align}\notag
	\lambda_{i_{1},i_{2},...,i_{d}}=\pi^{2}\sum_{k=1}^{d} i_{k}^{2}, \quad and~~ \varphi_{i_{1},i_{2},...,i_{d}}(\bm{x})=2\prod_{k=1}^{d}\sin (i_{k}\pi x_{k}),\quad i_{1},i_{2},...,i_{d} \in \mathbb{N}^{+}.
	\end{align}
	From (\ref{equ:4.1}), the exact solution $u(\bm{x})$ can be constructed by taking appropriate $f(\bm{x})$.
	\begin{example}\label{exa:4.1}(Smooth solution)~~%Example 1
		Taking $f(\bm{x})=\prod\limits_{k=1}^{d}\sin (n \pi x_{k})$ in equation (\ref{equ:2.1}) leads to
		\begin{align}\notag
			\widehat{f}_{i_{1},i_{2},...,i_{d}}=\frac{1}{2}\prod_{k=1}^{d}\left(\frac{\sin (n \pi-i_{k} \pi)}{n\pi-i_{k} \pi}-\frac{\sin (n\pi+i_{k} \pi)}{n\pi+i_{k} \pi}\right),\quad n \in \mathbb{N}^{+}.
		\end{align}
		It can be found that $\widehat{f}_{i_{1},i_{2},...,i_{d}}$ is $\frac{1}{2}$ in the sense of limitation if $i_{1}=i_{2}=\cdots=i_{d}=n$; otherwise, it will vanish. Hence we have from (\ref{equ:4.1}) that
		\begin{align}\notag
			u(\bm{x}) =\frac{1}{\left(dn^{2}\pi^{2}+\gamma\right)^{s}}\prod_{k=1}^{d}\sin (n\pi x_{k}), \quad n \in \mathbb{N}^{+}.
		\end{align}	
	\end{example}
	In Example \ref{exa:4.1}, we take $n=2$ and $d=3$. The error is measured by $\Vert e_{h}\Vert_{2}=\Vert u-u_{h}\Vert$ in $L^{2}$-norm, and the convergence rate is computed by $r_{2}=\log_{2}\left(\frac{\Vert e_{h}\Vert_{2}}{\Vert e_{h/2}\Vert_{2}}\right)$. Hence the expected accuracy is second order in linear FEM, and the expected accuracy is fourth order in CDM. The results are shown in Tables \ref{tab:4.1} and \ref{tab:4.2}.
	\begin{example}\label{exa:4.2}(Singular solution)~~%Example 2
		Let $f(\bm{x})=1$, $\gamma=1$ and $d=2$ in (\ref{equ:2.1}).
	\end{example}
		\begin{table*}[htb]\small
		\setlength{\abovecaptionskip}{0cm}
		\caption{\label{tab:4.3} $L^{2}$ errors and the convergence rates of Example \ref{exa:4.2} in FEM and CDM.}
		\centering
		\begin{tabular}{p{0.4cm} p{0.5cm} p{1.5cm} p{0.8cm}p{1.5cm}p{0.8cm} p{0.4cm} p{0.5cm} p{1.5cm} p{0.8cm}p{1.5cm}p{0.8cm}}
			\hline
			$s$ & $\mathcal{M}$  &  FEM         &$r_{2}$  &  CDM       &$r_{2}$&$s$   &$\mathcal{M}$   &  FEM      &$r_{2}$&  CDM       &$r_{2}$\\
			\hline 
	   0.5  & 64   &   3.996E-04  &         & 3.039E-04  &       & 0.9  & 64   & 1.984E-05 &       & 4.232E-05  &     \\
			& 128  &   1.377E-04  &  1.537  & 1.089E-04  & 1.480 &      & 128  & 4.935E-06 & 2.008 & 1.098E-05  &1.946\\
			& 256  &   4.798E-05  &  1.521  & 3.879E-05  & 1.490 &      & 256  & 1.231E-06 & 2.003 & 2.810E-06  &1.966\\
			& 512  &   1.682E-05  &  1.512  & 1.376E-05  & 1.495 &      & 512  & 3.075E-07 & 2.001 & 7.131E-07  &1.979\\
			& 1024 &   5.920E-06  &  1.507  & 4.873E-06  & 1.498 &      & 1024 & 7.684E-08 & 2.001 & 1.800E-07  &1.986\\
			& 2048 &   2.088E-06  &  1.504  & 1.725E-06  & 1.499 &      & 2048  & 1.921E-08 & 2.000 & 4.528E-08  &1.991\\
			&      &              &         &            &       &      &      &           &       &            &\\
			%\hline
	   1.3  & 64   &   3.649E-06  &         & 1.007E-05  &       & 1.7  & 64   & 5.738E-07 &       & 2.860E-06  &     \\
			& 128  &   9.124E-07  &  2.000  & 2.524E-06  & 1.997 &      & 128  & 1.436E-07 & 1.999 & 7.155E-07  &1.995\\
			& 256  &   2.281E-07  &  2.000  & 6.313E-07  & 1.999 &      & 256  & 3.590E-08 & 2.000 & 1.789E-07  &2.000\\
			& 512  &   5.702E-08  &  2.000  & 1.579E-07  & 2.000 &      & 512  & 8.976E-09 & 2.000 & 4.473E-08  &2.000\\
			& 1024 &   1.425E-08  &  2.000  & 3.947E-08  & 2.000 &      & 1024 & 2.244E-09 & 2.000 & 1.118E-08  &2.000\\
			& 2048 &   3.563E-09  &  2.000  & 9.867E-09  & 2.000 &      & 2048 & 5.610E-10 & 2.000 & 2.796E-09  &2.000\\
			\hline
		\end{tabular}
	\end{table*} 
	In Example \ref{exa:4.2}, the $L^{2}$ error can be measured by $\Vert e_{h}\Vert_{2}=\Vert u_{h}-u_{h/2}\Vert$, the convergence rate is given by $r_{2}=\log_{2}\left(\frac{\Vert e_{h}\Vert_{2}}{\Vert e_{h/2}\Vert_{2}}\right)$, and the expected accuracy is $\mathcal{O}(h^{2s+\frac{1}{2}})$ and $\mathcal{O}(h^{2})$ for $s<0.75$ and $s\geq 0.75$, respectively. The numerical results are presented in Table \ref{tab:4.3}, and the numerical accuracy is consistent with \cite{Dus2018AF}. It is noteworthy that the error is measured in $L^{2}$-norm, however, the numerical accuracy will be $\mathcal{O}(h^{2s})$ for $s\in \left(0,1\right)$ if using $L^{\infty}$-norm to measure the error.
\subsection{Time-dependent problems}
	In this subsection, we shall apply the present double fast algorithm to time-space fractional diffusion problems. To this end, we also take $Q=128$ and $\epsilon=10^{-16}$ in (\ref{equ:3.3}).                                                     
\begin{example}\label{exa:4.3}%Example 3
	Let $g$ be time-dependent smooth function, and assume that the exact solution of (\ref{equ:1.1})
	\begin{align}\notag
		u(\bm{x},t) =\frac{g(t)}{\left(dn^{2}\pi^{2}+\gamma\right)^{s}}\prod_{k=1}^{d}\sin (n\pi x_{k}), \quad n \in \mathbb{N}^{+},
	\end{align}
	is chosen so that
	\begin{align}\notag
		f(\bm{x},t)=\left(\frac{\partial_{t}^{\alpha}g(t)}{\left(dn^{2}\pi^{2}+\gamma\right)^{s}}+\kappa_{s}g(t)\right)\prod_{k=1}^{d}\sin (n\pi x_{k}), \quad n \in \mathbb{N}^{+}.
	\end{align}
\end{example}

	In Example \ref{exa:4.3}, we first take $d=2$, $\gamma=1$, $n=1$, $\omega=1$ and $\kappa_{s}=0.1$. In addition, the $L^{\infty}$ error is computed at $t_{N}=T$ and the spatial numerical method uses CDM. Then the spatial numerical results are presented in Table \ref{tab:4.4} with $g(t)=t$, $N=5000$ is fixed, and the expected accuracy is $\mathcal{O}(h^{4})$ in space. The temporal numerical results are shown in Table \ref{tab:4.5} with $g(t)=t^{1.5}$ and $\mathcal{M}=500$, and the expected accuracy is $\mathcal{O}(\Delta t^{2-\alpha})$ in time.
		\begin{table*}[htb]\small
		\setlength{\abovecaptionskip}{0cm}
		\caption{\label{tab:4.4} Spatial numerical errors and the convergence rates of Example \ref{exa:4.3} in $L^{\infty}$-norm.}
		\centering
		\renewcommand\arraystretch{1.1}
		\begin{tabular}{p{1.5cm} p{1cm} p{2cm} p{2cm}p{2cm} p{2cm} p{2cm}}
			\hline
			$(s,\alpha)$&                 $\mathcal{M}$              &   5          & 10        & 20        & 40          & 80 \\
			\hline 
			$(1.0,0.8)$  & $\Vert e_{h}\Vert_{\infty}$   &   1.650E-05  & 1.127E-06 & 7.021E-08 & 4.385E-09   & 2.740E-10            \\
			&  $~~r_{\infty}$               &              &  3.872    &  4.004    & 4.001       & 4.000  \\
			$(0.8,0.8)$  & $\Vert e_{h}\Vert_{\infty}$   &   1.747E-05  & 1.193E-06 & 7.437E-08 & 4.645E-09   & 2.902E-10  \\
			& $~~r_{\infty}$                &              &  3.872    &  4.004    & 4.001       & 4.000 \\
			\hline
			$(0.8,0.6)$  & $\Vert e_{h}\Vert_{\infty}$   &   1.861E-05  & 1.271E-06 & 7.920E-08 & 4.947E-09   & 3.090E-10  \\
			& $~~r_{\infty}$                &              &  3.872    & 4.004     & 4.001       & 4.001  \\
			$(0.6,0.6)$  & $\Vert e_{h}\Vert_{\infty}$   &   1.730E-05  & 1.182E-06 & 7.363E-08 & 4.598E-09   & 2.872E-10            \\
			& $~~r_{\infty}$                &              &  3.872    &  4.004    & 4.001       & 4.001  \\
			\hline
			$(0.6,0.4)$  & $\Vert e_{h}\Vert_{\infty}$   &   1.869E-05  & 1.277E-06 & 7.955E-08 & 4.966E-09   & 3.087E-10  \\
			& $~~r_{\infty}$                &              &  3.872    &  4.004    & 4.002       & 4.008 \\
			$(0.4,0.4)$  & $\Vert e_{h}\Vert_{\infty}$   &   1.459E-05  & 9.969E-07 & 6.212E-08 & 3.876E-09   & 2.386E-10  \\
			& $~~r_{\infty}$                &              &  3.872    & 4.004     & 4.002       & 4.022  \\
			\hline
		\end{tabular}
	\end{table*}
	\begin{table*}[htb]\small
		\setlength{\abovecaptionskip}{0cm}
		\caption{\label{tab:4.5} Temporal numerical errors and the convergence rates of Example \ref{exa:4.3} in $L^{\infty}$-norm.}
		\centering
		\renewcommand\arraystretch{1.1}
		\begin{tabular}{p{1.5cm} p{1cm} p{2cm} p{2cm}p{2cm} p{2cm} p{2cm}}
			\hline
			$(s,\alpha)$&               $N$                &   20         & 40        & 80        & 160         & 320 \\
			\hline 
			$(1.0,0.8)$  & $\Vert e_{h}\Vert_{\infty}$   &   1.963E-04  & 8.601E-05 & 3.761E-05 & 1.643E-05   & 7.168E-06            \\
			&  $~~r_{\infty}$               &              &  1.190    &  1.193    & 1.195       & 1.197  \\
			$(0.8,0.8)$  & $\Vert e_{h}\Vert_{\infty}$   &   5.746E-04  & 2.540E-04 & 1.117E-04 & 4.894E-05   & 2.140E-05  \\
			& $~~r_{\infty}$                &              &  1.178    &  1.186    & 1.190       & 1.193 \\
			\hline
			$(0.8,0.6)$  & $\Vert e_{h}\Vert_{\infty}$   &   2.013E-04  & 7.723E-05 & 2.951E-05 & 1.125E-05   & 4.280E-06  \\
			& $~~r_{\infty}$                &              &  1.382    & 1.388     & 1.392       & 1.394  \\
			$(0.6,0.6)$  & $\Vert e_{h}\Vert_{\infty}$   &   5.165E-04  & 1.993E-04 & 7.643E-05 & 2.920E-05   & 1.113E-05            \\
			& $~~r_{\infty}$                &              &  1.374    &  1.383    & 1.388       & 1.392  \\
			\hline
			$(0.6,0.4)$  & $\Vert e_{h}\Vert_{\infty}$   &   1.653E-04  & 5.583E-05 & 1.874E-05 & 6.258E-06   & 2.084E-06  \\
			& $~~r_{\infty}$                &              &  1.566    &  1.575    & 1.582       & 1.586 \\
			$(0.4,0.4)$  & $\Vert e_{h}\Vert_{\infty}$   &   3.802E-04  & 1.288E-04 & 4.333E-05 & 1.450E-05   & 4.833E-06  \\
			& $~~r_{\infty}$                &              &  1.561    & 1.572     & 1.579       & 1.585  \\
			\hline
		\end{tabular}
	\end{table*}
	
	To illustrate the efficiency of the present double fast algorithm for solving time-space fractional diffusion problem in Example \ref{exa:4.3}, we also provide a comparison of CPU time between fast L1 (fL1) scheme and direct L1 scheme. Here, the fast linear FEM (fFEM) and fast fourth-order CDM (fCDM) using discrete sine transform are carried out in space, respectively, $s=\alpha=0.4$ is fixed and $g(t)=t^{1.5}$ is chosen.
		\begin{table*}[htb]\small
		\setlength{\abovecaptionskip}{0cm}
		\caption{\label{tab:4.6} A comparison of CPU time (seconds) between direct L1 and fast L1 for Example \ref{exa:4.3} with $\mathcal{M}\approx\sqrt{N}$.}
		\centering
		\renewcommand\arraystretch{1.1}
		\begin{tabular}{p{1.5cm} p{1cm} p{2cm} p{2cm}p{2cm} p{2cm} p{2cm}}
			\hline
			Scheme &               $N$             &   1000          & 2000         & 4000         & 8000           & 16000 \\
			\hline 
			fFEM-L1& $\Vert e_{h}\Vert_{\infty}$   &   1.841156E-05  & 9.292752E-06 & 4.570950E-06 & 2.305831E-06   & 1.155627E-06  \\
			&        cpu (s)                &   2.286         &  10.47       &  60.76       & 422.4          & 2801 \\
			fFEM-fL1& $\Vert e_{h}\Vert_{\infty}$   &   1.841157E-05  & 9.292755E-06 & 4.570953E-06 & 2.305844E-06   & 1.155638E-06  \\
			&        cpu (s)                &   1.168         &  5.473       &  22.84       & 116.2          & 394.9 \\
			\hline
			fCDM-L1& $\Vert e_{h}\Vert_{\infty}$   &   7.974488E-07  & 2.641280E-07 & 8.715063E-08 & 2.878578E-08   & 9.505049E-09  \\
			&        cpu (s)                &   2.322         &  10.64       &  60.13       & 393.7          & 2660 \\
			fCDM-fL1& $\Vert e_{h}\Vert_{\infty}$   &   7.974468E-07  & 2.641250E-07 & 8.714802E-08 & 2.877249E-08   & 9.493764E-09  \\
			&        cpu (s)                &   1.182         &  5.357       &  21.58       & 112.0          & 449.9 \\
			\hline
		\end{tabular}
	\end{table*}
	
	\begin{table*}[htb]\small
		\setlength{\abovecaptionskip}{0cm}
		\caption{\label{tab:4.7} Optimal temporal numerical errors defined in Theorem \ref{the:3.4} and the convergence orders of Example \ref{exa:4.3}.}
		\centering
		\renewcommand\arraystretch{1.1}
		\begin{tabular}{p{1.5cm} p{1cm} p{2cm} p{2cm}p{2cm} p{2cm} p{2cm}}
			\hline
			$(\alpha,\omega)$&               $N$                &   20         & 40        & 80        & 160         & 320 \\
			\hline 
			$(0.9,1.222)$  & Error   &   8.574E-04  & 4.502E-04 & 2.3250E-04 & 1.185E-04   & 5.972E-05            \\
			&  Order               &              &  0.9292    &  0.9533    & 0.9726       & 0.9884  \\
			$(0.8,1.500)$  & Error   &   1.074E-03  & 5.156E-04 & 2.427E-04 & 1.126E-04   & 5.161E-05  \\
			& Order                &              &  1.059    &  1.087    & 1.108       & 1.125 \\
			\hline
			($0.7,1.857)$  & Error   &   1.025E-03  & 4.506E-04 & 1.940E-04 & 8.223E-05   & 3.447E-05  \\
			& Order                &              &  1.186    & 1.216     & 1.238       & 1.254  \\
			$(0.6,2.333)$  & Error   &   8.860E-04  & 3.583E-04 & 1.419E-04 & 5.539E-05   & 2.142E-05            \\
			& Order                &              &  1.306    &  1.336    & 1.357       & 1.371  \\
			\hline
			$(0.5,3.000)$  & Error   &   7.325E-04  & 2.743E-04 & 1.007E-04 & 3.655E-05   & 1.316E-05  \\
			& Order                &              &  1.417    &  1.445    & 1.463       & 1.473 \\
			$(0.4,4.000)$  & Error   &   5.932E-04  & 2.077E-04 & 7.158E-05 & 2.445E-05   & 8.302E-06  \\
			& Order                &              &  1.514    & 1.537     & 1.550       & 1.558  \\
			\hline
		\end{tabular}
	\end{table*}
	
	From Table \ref{tab:4.6}, we observed that the numerical accuracy of both L1 scheme and fL1 scheme are almost the same. For scheme using FEM, the error order exhibits $\mathcal{O}(h^{2}+N^{-1})$, which is due to the fact that the theoretical accuracy $\mathcal{O}(h^{2}+N^{-1.6})$ and $h=N^{-0.5}$ lead to that the spatial convergence quickly reaches optimal order but the order is merely 1 in time. While for scheme employing CDM, the error order enjoys $\mathcal{O}(h^{3.2}+N^{-1.6})$. This is because the temporal convergence achieves optimal order $2-\alpha$ early so that the order is not better than $4-2\alpha$ in space. For computational speed, the schemes using fL1 are significantly faster. These results demonstrate that the present double fast algorithm really can reduce the CPU time and save the memory.
	
	Table \ref{tab:4.7} presents a set of temporal numerical results which are mainly to illustrate the Theorem \ref{the:3.4}. The choices of all parameters are same as Table \ref{tab:4.6}. However, we shall consider a weak smooth time-dependent function $g(t)=t^{\alpha}$ in Example \ref{exa:4.3} satisfying time regularity Assumption \ref{ass:3.1} and also fix $\mathcal{M}=100$ to obtain the desired convergence orders. It can be noted that the choice of graded temporal mesh factor $\omega=(2-\alpha)/\alpha$ brought an optimal $\mathcal{O}(N^{-(2-\alpha)})$ convergence in time direction. Hence, these results are consistent with our theoretical analysis.
\begin{example}\label{exa:4.4}%Example 4
	In the example, we study the coarsening dynamics of the following time-space fractional phase-field problem
	\begin{align}\notag
	\left\{\begin{array}{l}
		\partial_{t}^{\alpha} u(\bm{x}, t)=-(-\Delta)^{\beta}\left(\varepsilon^{2s}(-\Delta+\gamma\mathbb{I})^{s}u(\bm{x},t)+F(u)\right), \quad \bm{x} \in \Omega, ~t>0, \\
		u(\bm{x}, t)=g(\bm{x}, t), \quad \bm{x} \in \partial\Omega, ~t \geq 0, \\
		u(\bm{x}, 0)=u_{0}(\bm{x}),\quad \bm{x} \in \Omega,
	\end{array}\right.
	\end{align}
	where $u$ is the scaled concentration subject to the given initial value $u_{0}$ and boundary data $g$, the nonlinear term $F(u)=u^{3}-u$, $0<\alpha, s,\beta<1$, $\gamma\geq 0$ and $\varepsilon>0$ is an interfacial parameter describing the thickness of the phase boundary.
\end{example}
	\begin{figure}[htb] %width=10cm, height=10cm
		\centering
		\includegraphics[width=0.85\textwidth]{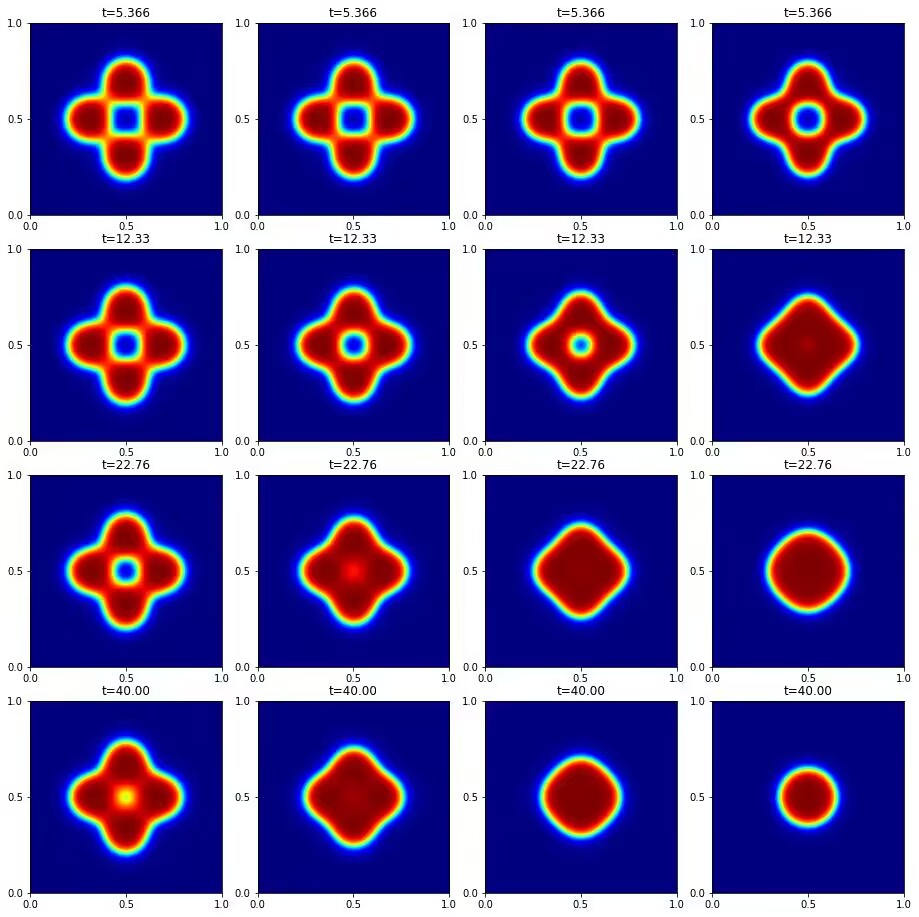} 
		\caption{Process of the evolution of time-space fractional Allen-Cahn equation with the kissing bubbles: $\alpha=0.4$, $0.6$, $0.8$, $1$ from left to right, respectively.}
		\label{img1}
	\end{figure}
	\begin{figure}[htb] %width=10cm, height=10cm
		\centering
		\includegraphics[width=0.85\textwidth]{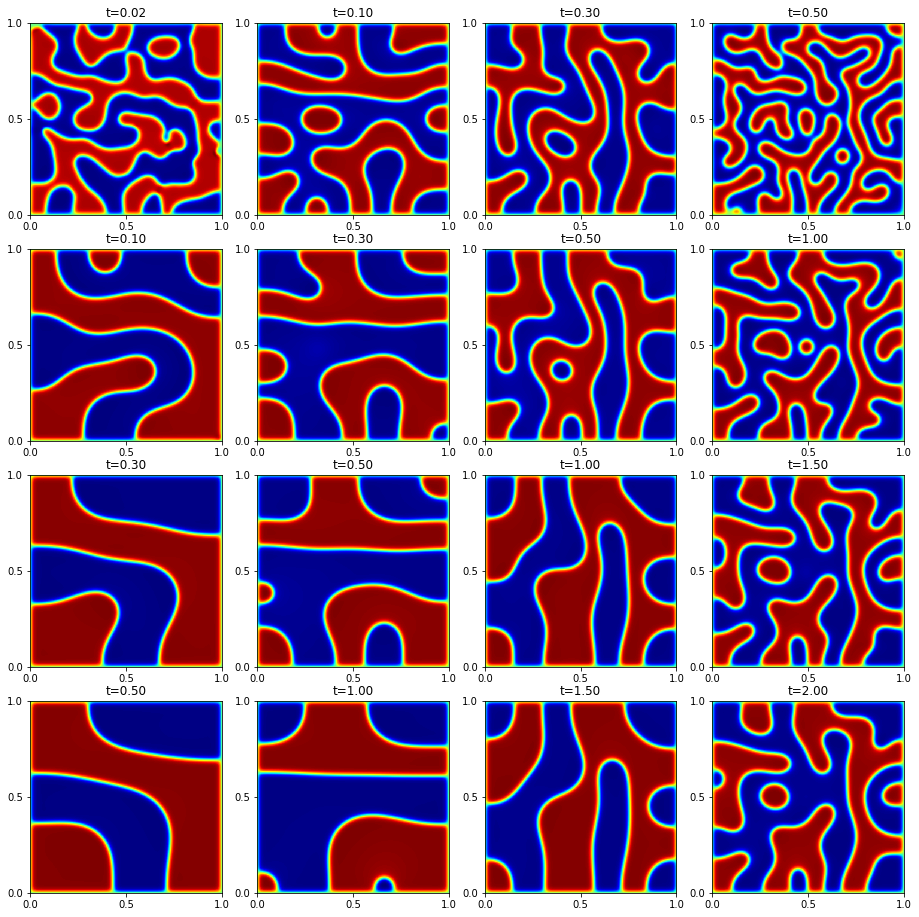}
		\caption{Snapshots of the numerical solutions of time-space fractional Cahn-Hilliard equation with $\beta=1.0$, $0.8$, $0.6$, $0.4$ from left to right, respectively.}
		\label{img2}
	\end{figure}
	
	In Example \ref{exa:4.4}, fast time-stepping L1 method is used, and linear FEM for spatial discretization is performed.
	
	We first consider $\beta=0$, $\gamma=1$, $\varepsilon=0.02$ and $g(\bm{x},t)=-1$ so that the Example \ref{exa:4.4} collapses to time-space fractional Allen-Cahn equation. Then we choose the initial value
	\begin{align}\notag
		u_{0}(\bm{x})=-\prod_{i=1}^{4}\tanh\left(\frac{|\bm{x}-\bm{x}_{i}|^{2}-0.015}{0.5\epsilon}\right),\quad \bm{x}\in \Omega,
	\end{align}
	where $\left(x_{1},y_{1}\right)=\left(0.7,0.5\right)$, $\left(x_{2},y_{2}\right)=\left(0.3,0.5\right)$, $\left(x_{3},y_{3}\right)=\left(0.5,0.7\right)$ and $\left(x_{4},y_{4}\right)=\left(0.5,0.3\right)$. In addition, take $u=\bar{u}-1$ such that the inhomogeneous Dirichlet boundary condition turns a homogeneous case.
	
	In the numerical test, to investigate the process of the evolution of four kissing bubbles for $\alpha=1, 0.8, 0.6, 0.4$ at different time level, we fix $s=0.8$, $h=2^{-10}$, $T=40$, $N=1600$ and $\omega=1.2$. From Figure \ref{img1}, we observed that four bubbles will shrink and turn into a single bubble with different rate of change, however, the process of the evolution of the bubble is slower for smaller $\alpha$.
	
	Next, we study the time-space fractional Cahn-Hilliard equation in Example \ref{exa:4.4} with $s=1$. Here, we choose initial value as a random perturbation uniformly distributed in $[-0.05, 0.05]^{2}$, and take $g(\bm{x},t)=0$, $\varepsilon=0.01$, $\omega=1$, $\alpha=0.8$, $\Delta t=0.001$ and $h=2^{-9}$. 
	
	We consider different $\beta$ to study the effect on coarsening behavior of the time-space fractional Cahn-Hilliard equation. Figure \ref{img2} shows
	the snapshots of the numerical solutions at different time for $\beta=1.0$, $0.8$, $0.6$ and $0.4$ from first column to last one, respectively, and the present numerical simulations use the same initial data. We found that the change of the states of solutions is significant, and the solutions would converge to different steady states for different values of $\beta$. 
	
	From Figures \ref{img1} and \ref{img2}, we observed that the proposed double fast algorithm can efficiently simulate the time-space fractional phase-field models. Some similar simulations can be seen in \cite{Bul2019SE}.
\section{Conclusion}
	In this paper, we proposed both commendably concise and implementation-friendly double fast algorithm for solving time-space fractional diffusion problems. We first employed linear finite element or fourth-order compact difference method combining with matrix transfer technique to approximate spectral fractional Laplacian, and then derived semi-discrete scheme of time-space fractional diffusion equation. In addition, a full discrete scheme using fast time-stepping L1 method on graded time mesh was provided, and then its $\alpha$-robust stability and convergence analyses were discussed. Numerical simulations showed that the present algorithm really could reduce computation cost and memory requirement. For instance, $\mathcal{M}=256$ and $\mathcal{M}=1024$ were taken in Examples \ref{exa:4.1} and \ref{exa:4.4}, the scales of algebraic matrix were over 16 million and one million, respectively. It is not realistic to obtain these results when using direct numerical methods.
	
	Our spatial algorithm is based on specific Cartesian meshes, and thus may not be directly extended to more complex regions. The design on the corresponding fast algorithms remains open and certainly should be a subject of future research. In the end, since the use of Lemma \ref{lem:3.6}, the error bound of Theorem \ref{the:3.4} exhibits an $\alpha$-nonrobust behavior, namely, the bound will blow up as $\alpha\to 1^{-}$.
\section*{Acknowledgments}
	The authors are very thankful to the editors and reviewers for carefully reading the paper, their suggestions and comments have improved the quality of the paper.
\section*{Data availability}
    Data will be made available on request.
\section*{Declarations}
    The authors declare that there is no conflict of interest.
\section*{Appendix}
	Here, we introduce a linear interpolation function of the solution $u(\zeta)$:
	\begin{align}\notag
		L_{k}(\zeta)=\frac{t_{k}-\zeta}{\Delta t}u(t_{k-1})+\frac{\zeta-t_{k-1}}{\Delta t}u(t_{k}),\quad \zeta\in (t_{k-1},t_{k}),\quad k=1,2,...,N,
	\end{align}
	where $\Delta t=t_{k}-t_{k-1}$. Define a second order difference quotient of $u(t)$ at $t_{k-1}$, $t_{k}$ and $\zeta$:
	\begin{align}\notag
		u[t_{k-1},t_{k},\zeta]=\frac{u[t_{k},\zeta]-u[t_{k-1},t_{k}]}{\zeta-t_{k-1}},\quad and~~ u[t_{i},t_{j}]=\frac{u(t_{j})-u(t_{i})}{t_{j}-t_{i}}.
	\end{align}
	According to Newton interpolation, we have
	\begin{align}\notag
		u(\zeta)=L_{k}(\zeta)+R_{k}(\zeta),\quad R_{k}(\zeta)=(\zeta-t_{k-1})(\zeta-t_{k})u[t_{k-1},t_{k},\zeta],\quad \zeta\in (t_{k-1},t_{k}),\quad k=1,2,...,N.
	\end{align}
	Then we derive direct L1 scheme as follows:
	\begin{align}\notag
		\partial_{t}^{\alpha}u(t_{n})=\sum_{k=1}^{n}c_{n,k}(u(t_{k})-u(t_{k-1}))+\bar{r}^{n},\quad n=1,2,...,N,
	\end{align}
	where $c_{n,k}$ and $\bar{r}_{n}$ are respectively given by
	\begin{align}\notag
		c_{n,k}=\frac{(t_{n}-t_{k-1})^{1-\alpha}-(t_{n}-t_{k})^{1-\alpha}}{\Delta t\Gamma(2-\alpha)},\quad \bar{r}_{n}=\frac{1}{\Gamma(1-\alpha)}\sum_{k=1}^{n}\int_{t_{k-1}}^{t_{k}}(t_{n}-\zeta)^{-\alpha}R_{k}^{'}(\zeta)d\zeta,\quad n=1,2,...,N.
	\end{align}
	
	Now, we shall conclude a result as follows.
	\begin{proposition}\label{pro:5.1} Let $m=0,1$, and $0\leq\delta\leq 1$. For $u(t)\in C^{m,\delta}([0,T])$, we have
		\begin{align}\notag
			\max_{1\leq n\leq N}\left|\bar{r}_n\right|\leq C_{m}\Delta t^{m+\delta-\alpha}\|u\|_{C^{m,\delta}([0, T])},
		\end{align}
		where if $m=0$, then $\delta\in (\alpha,1]$ and $C_{0}=\frac{\alpha}{\Gamma(1-\alpha)}\left(\frac{2}{\alpha}+\frac{1}{\delta-\alpha}+\frac{1}{1-\alpha}\right)$; and if $m=1$, then $\delta\in \left[0,1\right]$ and $C_{1}=\frac{1}{\Gamma(2-\alpha)}$. In addition, $C^{m, \delta}([0, T])$ denotes Hölder continuous function space equipped with norm
		\begin{align}\notag
			\|u\|_{C^{m,\delta}([0, T])}=\sum_{i=0}^{m}\|\partial_{t}^{i}u\|_{C([0, T])}+\sup _{t_1, t_2 \in(0, T); t_1 \neq t_2} \frac{\left|\partial_{t}^{m}u\left(t_1\right)-\partial_{t}^{m}u\left(t_2\right)\right|}{\left|t_1-t_2\right|^\delta}, \quad 0\leq \delta\leq 1.
		\end{align}
	\end{proposition}
	\begin{proof} For $m=1$, the limited smoothness result has been proved in \cite[Theorem 2.1]{Zht2020TH} and it can be noted that the solution $u(t)$ has at least the $C^1{}$ regularity in time. Next, we consider weak smooth case $m=0$. Using the following relation
		\begin{align}
			R_k(\zeta)=\left(\zeta-t_k\right)\left(u\left[t_k,\zeta\right]-u\left[t_{k-1}, t_k\right]\right),\quad \zeta \in\left(t_{k-1}, t_k\right),\quad R_k\left(t_{k-1}\right)=R_k\left(t_k\right)=0,\notag
		\end{align}
		and integration by parts, we have
		\begin{align}
			\left|\bar{r}_n\right|=\frac{1}{\Gamma(1-\alpha)}\left|\sum_{k=1}^{n-1}-\alpha \int_{t_{k-1}}^{t_k}\left(t_n-\zeta\right)^{-\alpha-1} R_k(\zeta) d \zeta+\int_{t_{n-1}}^{t_n}\left(t_n-\zeta\right)^{-\alpha} R_n^{\prime}(\zeta)d\zeta\right|.\notag
		\end{align}
		Now we shall use integration by parts again to obtain
		\begin{align}
			\int_{t_{n-1}}^{t_n}\left(t_n-\zeta\right)^{-\alpha} R_n^{\prime}(\zeta)d\zeta=-\alpha\int_{t_{n-1}}^{t_n}\left(t_n-\zeta\right)^{-\alpha-1} R_n(\zeta) d\zeta.\notag
		\end{align}
		Hence, we can bound
		\begin{align}
			\left|\bar{r}_n\right| \leq \frac{\alpha}{\Gamma(1-\alpha)}\left(\sum_{k=1}^{n-1} \int_{t_{k-1}}^{t_k}\left(t_n-\zeta\right)^{-\alpha-1}\left|R_k(\zeta)\right|d \zeta+\int_{t_{n-1}}^{t_n}\left(t_n-\zeta\right)^{-\alpha-1}\left|R_n(\zeta)\right|d\zeta\right).\notag
		\end{align}
		A direct triangle inequality and the definition of Hölder norm in $C^{0, \delta}([0, T])$ space give
		\begin{align}
			\left|R_k(\zeta)\right|&=\left|\left(\zeta-t_k\right)\left(u[t_{k},\zeta]-u[t_{k-1},t_{k}]\right)\right|\leq
			\left(t_k-\zeta\right)\left(\left|u[t_{k},\zeta]\right|+\left|u[t_{k-1},t_{k}]\right|\right)\notag\\
			&=\left(t_k-\zeta\right)\left(\left|\frac{u(\zeta)-u(t_{k})}{\zeta-t_{k}}\right|+\left|\frac{u(t_{k})-u(t_{k-1})}{t_{k}-t_{k-1}}\right|\right)\notag\\
			&\leq\left(\left(t_{k}-\zeta\right)^{\delta}+\Delta t^{\delta-1}\left(t_k-\zeta\right)\right)
			\|u\|_{C^{0, \delta}([0, T])},\quad \zeta \in\left(t_{k-1}, t_k\right).\notag
		\end{align}
		Combining $\left|\bar{r}_n\right|$ and $\left|R_k(\zeta)\right|$ leads to
		\begin{align}
			\left|\bar{r}_n\right| & \leq \frac{2\Delta t^{\delta}\alpha}{\Gamma(1-\alpha)}\int_0^{t_{n-1}} \left(t_n-\zeta\right)^{-\alpha-1} d \zeta
			\|u\|_{C^{0, \delta}([0, T])} \notag\\
			&+\frac{\alpha}{\Gamma(1-\alpha)}\int_{t_{n-1}}^{t_n}\left(t_n-\zeta\right)^{-\alpha-1}\left(\left(t_{n}-\zeta\right)^{\delta}+\Delta t^{\delta-1}\left(t_n-\zeta\right)d\zeta\right)\|u\|_{C^{0, \delta}([0, T])}\notag\\
			&=\frac{2\Delta t^{\delta}\alpha}{\Gamma(1-\alpha)}\left(\frac{\Delta t^{-\alpha}-t_{n}^{-\alpha}}{\alpha}\right)\|u\|_{C^{0, \delta}([0, T])}
			+\frac{\alpha}{\Gamma(1-\alpha)}\left(\frac{\Delta t^{\delta-\alpha}}{\delta-\alpha}+\Delta t^{\delta-1}\frac{\Delta t^{1-\alpha}}{1-\alpha}\right)\|u\|_{C^{0, \delta}([0, T])}\notag\\
			& \leq \frac{\alpha \Delta t^{\delta-\alpha}}{\Gamma(1-\alpha)}\left(\frac{2}{\alpha}+\frac{1}{\delta-\alpha}+\frac{1}{1-\alpha}\right)\|u\|_{C^{0,\delta}([0, T])},\quad \delta>\alpha.\notag
		\end{align}
		Hence, we end this proof.
	\end{proof}
	
	We can note from Proposition \ref{pro:5.1} that the solution $u(t)$ requires at least the $C^{0,\alpha+\epsilon}$ regularity with arbitrary small $\epsilon>0$ in time when $m=0$. The result is consistent with our normal cognition. However, for properly smooth data $f$ and $u_{0}$, though the solution of time fractional diffusion equation is sufficiently smooth for $t > 0$, it usually exhibits a weakly singular behavior at $t = 0$. For more details, see \cite{Zht2020TH,Stm2017EA,Chh2021BU}.
\bibliographystyle{unsrt}
\bibliography{refs}
\end{document}